 \documentclass[10pt,twocolumn,twoside]{IEEEtran} 

\IEEEoverridecommandlockouts                              

\usepackage[utf8]{inputenc}
\usepackage{amsmath,amssymb,bbm, amsfonts, bm}
\usepackage{graphicx, pgf, color, subfigure, scalefnt}
\PassOptionsToPackage{hyphens}{url}\usepackage{hyperref}
\usepackage{ enumerate, cite}
\usepackage{array}
\usepackage{booktabs}
\usepackage{multirow}
\usepackage{adjustbox}

\makeatletter
\renewcommand*{\@opargbegintheorem}[3]{\trivlist
      \item[\hskip \labelsep{\itshape $\quad$#1\ #2}] \emph{(#3):}\ }
\makeatother

\def\P{\mathbbm{P}}
\def\E{\mathbbm{E}}

\def\reals{\mathbb{R}}

\newtheorem{theorem}{Theorem}
\newtheorem{lemma}{Lemma}
\newtheorem{proposition}{Proposition}

\newtheorem{remark}{Remark}
\newtheorem{definition}{Definition}

\def\Qcal{\mathcal Q}
\def\q2x{\sigma}

\def\Pcal{\mathcal P}
\def\fb{\bar f}
\def\ones{\mathbf{1}}
\def\Prob{\mathbb{P}}
\def\Ical{\mathcal I}
\def\IDA{\Ical^\mathrm{DA}}
\def\IRT{\Ical^\mathrm{RT}}
\def\Lcal{\mathcal L}

\def\Zbb{\mathbb{Z}}
\def\Qbb{\mathbb{Q}}

\def\Ecal{\mathcal E}
\def\Vcal{\mathcal V}

\usepackage{url}
\usepackage{braket}
\newcommand{\rev}[1]{\textcolor{black}{#1}}
\newcommand{\rv}[1]{{#1}}

\title{
Flexible Market for Smart Grid:\\
Coordinated Trading of Contingent Contracts
}

\author{Junjie Qin,~\IEEEmembership{Student Member, IEEE,}
Ram Rajagopal,~\IEEEmembership{Member, IEEE,}
and Pravin Varaiya,~\IEEEmembership{Life Fellow, IEEE}
\thanks{J. Qin is with Institute for Computational and Mathematical Engineering, Stanford University, Stanford, CA, 94305. Email:
        {\tt\small jqin@stanford.edu}}%
\thanks{R. Rajagopal is with the Department of Civil and Environmental Engineering, Stanford University, Stanford, CA, 94305.
        Email: {\tt\small ramr@stanford.edu}}%
\thanks{P. Varaiya is with the Department of Electrical Engineering and Computer Science, U.C. Berkeley, Berkeley, CA, 94720.
        Email: {\tt\small varaiya@berkeley.edu}}%
\thanks{This work was supported in part by  the National Science Foundation under  Awards ECCS-1549945 and  ECCS-1554178 CAREER, in part by ARPA-e under Award No. DE-AR0000697, 
and in part by the Satre Family Fellowship.}
}

\begin{document}


\maketitle
\thispagestyle{empty}
\pagestyle{empty}

\begin{abstract}

A coordinated trading process is proposed as a design for an electricity market with significant uncertainty, perhaps from renewables.  In this process\rv{,} groups of agents propose to the system operator (SO) a contingent buy and sell trade that is balanced, i.e. the sum of demand bids and the sum of supply bids are equal.  The SO accepts the proposed trade if no network constraint is violated or curtails it until no violation occurs.  Each proposed trade is accepted or curtailed as it is presented.  The SO also provides guidance to help future proposed trades meet network constraints. 
The SO does not set prices, and there is no requirement that different trades occur simultaneously or clear at uniform prices.  Indeed, there is no price-setting mechanism.  However, if participants exploit opportunities for gain, the trading process will lead to an efficient allocation of energy and to the discovery of locational marginal prices (LMPs).  The great flexibility in the proposed trading process and the low communication and control burden on the SO may make the process suitable for coordinating producers and consumers in the distribution system.

\end{abstract}

\section{Introduction}\label{sec:intro}

The electric power grid increasingly relies  on distributed and variable energy sources.  Integration of these new  sources is helped by a market that facilitates matching intermittent supply and flexible demand  \cite{doe2030, ren212014}. Today the system operator (SO) achieves resource adequacy, congestion management and efficiency through reserve requirements, day-ahead (DA) and real-time (RT) markets, and centralized dispatch of standard energy commodities, namely specified amount of energy delivered at specified nodes
at fixed prices  \cite{StoftBook, KirschenBook}. The needs of  important participants cannot be adequately expressed in terms of these 
standard  commodities, and so the
 SO allows bilateral contracts (e.g. Google \cite{google}, GM \cite{gm} and Amazon \cite{amazon}), with contractual arrangements that are not known to the SO.  Over time, the rigidity of the standard commodity was more broadly felt and fitful accommodations were made by introducing new commodities, such as demand response, ramping, and capacity.  But given the legacy of the standard market this slow expansion of the SO's responsibility  cannot unlock the full contractual flexibility that participants may wish. In particular, it is challenging to repurpose today's market design to serve the needs of  distribution system operators (DSOs) who must coordinate participants with small distributed generation and controllable demand side devices, and who would benefit from differentiated micro-contracts (i.e. contracts whose volumes are of the order of kWhs). Possible examples of such differentiated contracts are (i) contracts for flexible amount of energy contingent on the realization of uncertain supply or demand, (ii) contracts to serve deferrable loads that consume a fixed amount of energy for (say) one hour but which could be scheduled for any hour of the day \cite{BitarLow,DDServices}, (iii) contracts that  favorably price  generation sources that are green or more flexible, and (iv) contracts that encourage local sharing among prosumers with solar PV and storage devices. Incorporating these differentiated contracts requires a significant deviation from an electricity market with a small number of standard commodities. 

In this paper we propose a more flexible alternative to the current market design, called \emph{coordinated multilateral trading}. In this design, participants trade among themselves according to  terms and conditions fashioned to suit their own purposes  like in today's over the counter  markets (OTC),
in contrast with exchanges for trading standardized commodities at transparent prices.  These are  \emph{contingent trades} as the amount of energy delivered is contingent on events or conditions specified in the contracts. Since the trades induce power flow, they must be coordinated to ensure that  network constraints are not violated.  The SO accomplishes this coordination task by curtailing trades if network constraints are violated, and publishing  information about the network state to guide participants regarding how subsequent trades can avoid overloading congested lines. Thus the proposed market design permits flexible contracts by allowing contingent trades while the SO maintains power system security. In today's design the SO computes an efficient dispatch  that respects network line constraints, but in the proposed design the SO is only concerned with  reliability, and the determination of an \emph{efficient} dispatch is left to self-interested  participants.

\subsection{Contributions and Organization}


Coordinated  trading of contingent contracts (described in Section~\ref{sec:trading}) is proposed as a flexible market mechanism in the context of electric power transmission system operation. 
We establish that the    trading process is well-defined and during each step of the process,  power system reliability is guaranteed though  the role of the SO is greatly simplified. Furthermore, we show that the trading process converges to an efficient dispatch, which meets a benchmark defined using social welfare maximization as in the centralized stochastic economic dispatch (Section~\ref{sec:eff}). We also show that this trading process discovers the optimal \emph{locational marginal prices} through the marginal costs of local participants (Section~\ref{sec:price}). Finally, we prove that the dispatch and prices identified from the trading process support an Arrow-Debreu equilibrium, a notion of competitive equilibrium under uncertainty (Section~\ref{sec:ad}). The trading process is illustrated with a simple two-bus example in Section~\ref{sec:eg}.  

\subsection{Related Literature}


In  studies of the standard electricity market  
the basic framework is a one- or two-settlement market (Day Ahead Market or DAM and Real Time Market or RTM) in a deterministic setting \cite{ Hogan1992, Baldick2004, 1033699, 5717327}.  In this framework, generators and load serving entities present supply and demand function bids to the system operator (SO); the SO then calculates the equilibrium as the generation and load schedule that maximizes social welfare (producer plus consumer surplus), subject to the constraint that flows on transmission lines are limited by their rated capacities.  This centralized calculation has the form of a mathematical programming problem   called the optimum power flow problem.  The dual variables at the optimum solution are called locational or nodal marginal prices or LMPs.  The LMP at a node is the marginal  cost of delivering  additional power at that node.  In a two-settlement market there are day ahead and real time LMPs.

In a stochastic context,  uncertainty is modeled by a probability distribution over a  set of scenarios.  Each scenario  has its specific supply and demand functions, and the SO finds the schedule that maximizes expected social welfare.  This schedule is contingent, since there is a different  schedule for each scenario \cite{THEC:THEC11, 7403136,doi:10.1287/opre.1090.0800, 6146389, 7447005, 7436360}. The complexity of the stochastic problem grows in three ways with the number of scenarios.  First, each demand and supply  bid  now is a function of prices \textit{and} scenarios, so the number of decision variables and LMPs will be multiplied by the number of scenarios,       thereby increasing  the SO's communication and computational burden. Second, there must be agreement among all participants about the probability distribution over the scenarios, which precludes heterogeneous beliefs or private information that can affect beliefs.  Third, participants must work out in advance the bids they will offer for each scenario and price vector.   This complexity has precluded real world implementation of the optimal stochastic power flow problem.  In the absence of contingent (stochastic) bids that permit risk mitigation and reduce volatility, stochastic perturbations in demand and supply may lead to the large variations in LMPs that are observed.

Two studies propose decentralized trading processes to replace the SO's centralized calculation. 
In \cite{Chao1996}, 
transmission rights are privately owned; the SO specifies ``marginal loading factors'' that is, the amount of  capacity on every transmission line that must be purchased by every proposed bilateral transaction. Transmission  prices are adjusted iteratively in steps as follows.  At any step nodal price differences  adjust to eliminate arbitrage profits from purchasing energy at one node and selling at another.  Given nodal prices, transmission  prices then are adjusted to increase rents, subject to the competitiveness condition that the transmission price for a line with excess capacity must be zero.  The iteration converges in the limit to the welfare maximizing solution, and the nodal prices converge to the LMPs.  

Our proposed design is closer to the decentralized multilateral trading process in \cite{Wu199975} \rv{and generalizes their trading process developed for single period deterministic electricity market into a setting with two periods and with uncertainty explicitly considered. 
}  In the multilateral trading process, groups of buyers and sellers propose to the SO a balanced trade, i.e. sum of buy bids equals sum of sell bids.  The SO accepts the trade if (together with previously accepted trades) no transmission line constraint is violated.  Otherwise, the SO curtails the proposed trade until the violation is eliminated.  No price is announced.  The understanding is that the private terms and conditions of a trade (including monetary payments) are acceptable to all parties.  As in  \cite{Chao1996} 
the SO announces loading vectors to guide participants towards trades that do not violate line constraints.  It is shown that in case generators are motivated by profit maximization and buyers by utility maximization, the process will converge to
a social welfare maximum.

Two important distinctions between these  decentralized processes are worth noting.  First, in the language of mathematical programming,  \cite{Chao1996} describes a dual method, whereas \cite{ Wu199975} gives a primal method.
It is possible that at each step the iteration in \cite{Chao1996} is infeasible except in the limit, whereas each step in 
\cite{ Wu199975} is feasible and the process may be stopped at any point.  Second, even though it is decentralized, the
process in \cite{Chao1996} is synchronized: trades in each step must occur at the same time; but in \cite{ Wu199975} trades are asynchronous.

\subsection{Notation}
For a natural number $N$,  $[N]$  denotes the set $\{1,\dots, N\}$. Let $x \in \reals^{I\times J}$ be a matrix, with entries denoted by $x_{i,j}$, $i \in [I]$, $j\in [J]$. We use $x_i \in \reals^{J}$ to denote the vector $(x_{i,j})_{j\in [J]}$ and $x_j \in \reals^{I}$ to denote the vector $(x_{i,j})_{i \in [I]}$. For an Euclidean space $\reals^d$, we use $\ones \in \reals^d$ to denote the all-one vector. 

\section{Formulation}\label{sec:form}
\subsection{Network Model}
Consider a power network with $N$ buses and $L$ power lines with capacity constraints. 
\rv{The physics of power network is  described by the AC power flow equations, which is a set of quadratic equations relating the nodal complex voltages with nodal complex power injections \cite{wood2006power}. In this paper, we focus on the real power component of the complex power and adopt  a linearized DC  approximation to the AC power flow equations\footnote{\rv{The DC approximation does not incorporate voltage constraints which are important for distribution systems. A different linearization of the AC power flow, the linearized DistFlow model, should be used for distribution system applications \cite{19266, 25627}. }}\cite{4956966}, as commonly used in current electricity markets\footnote{\rv{Some of today's transmission system operators  utilize a mixture of AC and DC model for their operation (e.g. ISO New England uses a DC optimal power flow with AC feasibility and PJM employs a DC-AC iteration) \cite{federal2011recent}. Although this paper assumes a linearized power flow model, many of our results can be extended to the AC power flow model with some work (cf. Section 7.6 and 7.7 of \cite{Wu199975}). }} \cite{1198270, 6935037, rau2003optimization}. 
Under such a linearization, the region of feasible nodal power injections has the form 
\begin{equation}
	\Pcal := \left\{ p\in \reals^{N}: -\bar r\le \widehat{H} p \le \bar r, \quad \ones^\top p = 0\right\}, 
\end{equation}
 where $\widehat{H} \in \reals^{L \times N}$ is the linear sensitivity matrix relating the line flows to the nodal injections, and $\bar r \in \reals^{L}$ records the capacity values of the  lines. To simplify the notation, we denote 
 \begin{equation}
 	H = \begin{bmatrix}
 \widehat H\\
 -\widehat H	
 \end{bmatrix},\quad 
\mbox{and} \quad \fb = \begin{bmatrix}
	\bar r\\
	\bar r
\end{bmatrix}, 
 \end{equation}
 so the region of feasible nodal injections  becomes
\begin{equation}\label{eq:Pcal}
\Pcal := \left\{ p\in \reals^{N}: Hp \le \fb, \quad \ones^\top p = 0\right\}.
\end{equation}
}The first inequality  in \eqref{eq:Pcal} models the line capacity constraints, while the second equality enforces power balance over the entire network. We will denote the rows of $H$, referred to as \emph{loading vectors} or shift factors, by $h_\ell^\top \in \reals^{1\times N}$. 
Throughout the paper, we assume that $\Pcal$ has a non-empty interior. 

\subsection{Uncertainty Model}
We consider the operation of the electricity market over two time periods,  the DA market and the RT market. 

We explicitly model the RT uncertainty as a finite collection of $S$ system scenarios, so  each scenario is indexed by $s\in [S]$ with probability  $\Prob(s)>0$.
We assume that the set of scenarios and \rv{the probabilities} are known to all  market participants and the SO and the realization of a scenario is publicly verifiable by all of them.
We could have the set of feasible injections $\Pcal$ depending on the scenario as well to model
transmission line failures, in which case $H$ and $\fb$ in 
\eqref{eq:Pcal} will be indexed by scenario $s$.  We do not do this to simplify the notation.

\subsection{Participant Model}\label{sec:part}
On each bus of the network $n \in [N]$, there resides a collection of electricity market participants denoted by $\Ical_n$, each of which is either an electricity producer or an electricity load. We model each market participant by her \rv{(or his)} RT \emph{power injection plan}, denoted by $p_i = (p_{i,s})_{s\in [S]}$, her local feasible power injection sets, denoted by $\Pcal_{i,s}$ such that $p_{i,s} \in \Pcal_{i,s}$ for all $i\in \Ical_n$ and $s\in [S]$,  and her von Neumann–Morgenstern utility function over such a plan, denoted by $U_i(p_i)$ and taking the form of
\begin{equation}\label{eq:vnmU}
U_i(p_i) =\rev{ \E [u_{i,s} (p_{i,s})] = \sum_{s\in [S]} \Prob(s) u_{i,s}(p_{i,s}),}
\end{equation}
where $u_{i,s} (p_{i,s})$ is the actual utility given scenario $s$. 
Throughout the paper, we assume a \emph{quasi-linear environment}, so that the utility function is linear in the amount of monetary payment of each market participant, i.e., 
\begin{equation}\label{eq:uis}
u_{i,s}(p_{i,s}) = m_{i,s}(p_{i,s}) + \tilde u_{i,s}(p_{i,s}),
\end{equation}
where $m_{i,s}(p_{i,s}) $ is the payment received by the participant in scenario $s$ and  $\tilde u_{i,s}(p_{i,s})$ is the utility associated with power injection $p_{i,s}$ as discussed in detail below.  We will assume that the utility function $\tilde u_{i,s}(\cdot)$ is concave for each $i\in \Ical$ and $s\in [S]$. 

For an electricity producer, the power injection is induced by the producer's possibly scenario-dependent electricity production so that $p_i \in \reals^S_+$. The feasible power injection sets model the generation limits, which could be scenario dependent in the case of renewable generation. Thus, we have $\Pcal_{i,s} = [0, \bar p_{i,s}] $, where $\bar p_{i,s}$ is the maximum possible power output in scenario $s$. The utility function is as defined in \eqref{eq:vnmU} and \eqref{eq:uis}, with 
\begin{equation}
 \tilde u_{i,s}(p_{i,s}) = - c_{i,s}(p_{i,s}),
\end{equation}
where $c_{i,s}(\cdot)$ is the cost function of  the generation plant. 

For an electricity load, the power injection is induced by the possibly scenario-dependent electricity consumption so that $p_i \in \reals^S_-$. Symmetrically with the producer, we have $\Pcal_{i,s} = [-\bar p_{i,s}, 0]$, where $\bar p_{i,s}$ is the maximum possible power demand in scenario $s$. The utility function is taken to be 
\begin{equation}
	 \tilde u_{i,s}(p_{i,s}) = b_{i,s}(p_{i,s}),
\end{equation}
where $b_{i,s}(\cdot)$ characterizes the benefit of using power by the particular load. \rv{For large loads (e.g. resellers), the benefit corresponds to the profit made from the given power consumption; for small loads such as individual consumers, the benefit function reflects the monetary value of consuming electricity and is a widely used device for modeling how power consumption varies with prices \cite{7858963, schweppe2013spot, 6266724}. Allowing the benefit function to be scenario dependent is useful for modeling e.g. demand response resources whose availability is not known a priori.}

 We partition the set $\Ical_n$ as $\Ical_n  = \IDA_n \cup \IRT_n$ such that $\IDA_n \cap \IRT_n = \emptyset$ and denote $\IDA = \cup_{n\in [N]} \IDA_n$ and $\IRT = \cup_{n\in[N]} \IRT_n$, where $ \IDA_n$ contains producers/loads connected to bus $n$ whose power injection has to be fixed in DA and cannot adapt to RT scenarios and $\IRT_n$ are those that can adapt to RT scenarios. We refer \rv{to} participants in $\IDA$ as \emph{DA participants} and those in $\IRT$  as \emph{RT participants}. 
Technically, the power injection plan of DA participants must satisfy the \emph{non-anticipation constraint}
 \begin{equation}
 	 p_{i,s} = p_{i,\tilde s}, \quad \mbox{for all } s, \tilde s \in [S]. 
 \end{equation}
 To \rv{simplify} the notation, let 
 \begin{equation}
 \bar \Pcal_i = \{p_i \in \reals^{S}: p_{i,s} \in \Pcal_{i,s}, s \in [S]\},
 \end{equation}
 and 
 \begin{equation}
 \Pcal_ i = \begin{cases}
 	\bar \Pcal_i \cap \{p_i \in \reals^S: p_{i,s} = p_{i,\tilde s },\,\, s,\tilde s\in [S]\}, &  i \in \IDA, \\ 
 	\bar \Pcal_i, & i \in \IRT.
 \end{cases}
 \end{equation}
Examples of DA participants include power plants that cannot ramp up or down following the RT uncertainty, such as coal-based generation plants, and loads that contract a fixed amount of consumption in each hour in DA. RT participants can either be variable generation sources or demand modeling e.g. renewable generation or random power consumption, or controllable generation or demand that can adapt to RT scenarios, such as fast-ramping gas generation or demand response resources. 

\subsection{Efficiency Benchmark}
A commonly used criteria for economic efficiency is \emph{Pareto optimality}. In a quasi-linear environment, \rev{it is equivalent to the following stochastic \emph{social welfare maximization} problem}
\begin{subequations}\label{opt:swm}
\begin{align}
	\mbox{maximize}\quad & U(p):= \sum_{i \in \Ical} U_i(p_i) \\
	\mbox{subject to}\quad & p_{i} \in \Pcal_{i}, \quad i \in \Ical, \label{opt:swm:local}\\
											& x_{n,s} = \sum_{i \in \Ical_n} p_{i,s}, \quad n\in [N],\,\, s \in [S], \label{opt:swm:pb}\\
											& x_s \in \Pcal, \quad s \in [S]. 
\end{align}
\end{subequations}
Notice that when the system does not take money from outside sources, we must have \emph{ex post budget adequacy}:
\begin{equation}
\sum_{i \in \Ical} m_{i,s} (p_{i,s}) \le 0.
\end{equation}
If  \emph{ex post budget balance} holds, i.e., $\sum_{i \in \Ical} m_{i,s} (p_{i,s}) = 0$, then the social welfare maximization program~\eqref{opt:swm} is equivalent to the \emph{stochastic economic dispatch} problem with the objective of~\eqref{opt:swm} replaced by
\begin{equation}
\sum_{i \in \Ical} U_i(p_i) = \rev{\sum_{i \in \Ical}  \E\left[ \tilde u_{i,s}(p_{i,s})\right]},
\end{equation}
where the summation  is the net sum of 
ex ante generation costs and load benefits as discussed in the previous subsection.


\section{Trading Process} \label{sec:trading}
The simplest market mechanism is one based on meeting and trading among self-interested agents. The electricity market is different in that centralized coordination has been commonly considered  essential to  ensure power system reliability constraints \eqref{eq:Pcal}. Indeed, completely decentralized trading without coordination could lead to line flows that violate their capacity limits and compromise the reliability of the system. As such, the standard power system market designs rely on a centralized clearing house (or market maker), referred to as system operator (SO), to solve an economic dispatch optimization in order to determine the generator schedules and electricity prices. When uncertainty from renewable generation is considered, the resulting stochastic economic dispatch problem is computationally more complex  and leads to increased communication requirement between the SO and market participants.  

\begin{figure}[htbp]
\centering
\includegraphics[scale = 0.4]{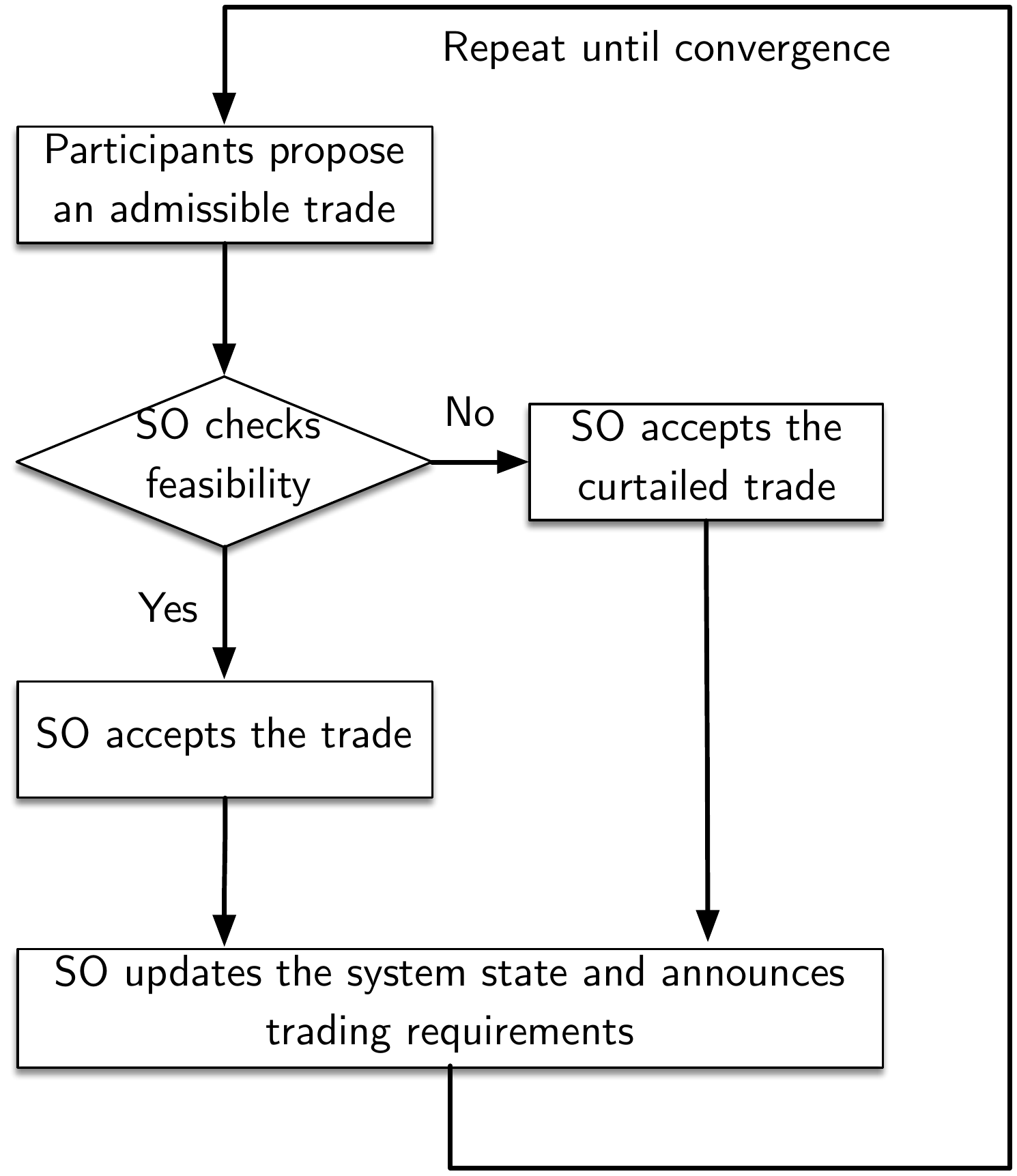}
\caption{Conceptual diagram for the trading process.}
\label{fig:flowchart}
\end{figure}

Wu and Varaiya \cite{Wu199975} propose a remarkably simple fix to make the free-market style meet-and-trade procedure respect the power system reliability constraints~\eqref{eq:Pcal}. The idea is to inject minimal amount of coordination, implemented by the SO, into the free trades so that the reliability (or feasibility) is guaranteed in every step of the trading process as shown in Figure~\ref{fig:flowchart}.  They also establish that the trading process achieves  economic efficiency in the limit. 
We will  generalize their  \emph{coordinated trading framework}  developed for single period deterministic electricity market into a setting with two periods and with uncertainty explicitly considered. 
\rv{Although we consider only a two-period market (consisting of a forward market, i.e. DA, and a delivery period, i.e. RT) below, the analysis readily extends to settings with multiple delivery periods. }

\rv{In this paper, we consider a setting where \emph{all} market participants trade exclusively in DA. This means that during the DA market, each DA participant $i \in \IDA$ trades and determines her power injection while each RT participant $i \in \IRT$ trades and determines her \emph{contingent power injection plan}. Since under the current setting with a \emph{complete} DA forward market, there is no need for RT re-trading, i.e., additional RT trading can not improve  social welfare, we assume that there is no trading  in real time.} 

\rv{We start with definitions for the trading process. Given the abstract nature of some of the definitions, examples demonstrating them are provided in Section~\ref{sec:eg} and linked here in footnotes.}

The premise of our trading system is that self-interested market participants will meet and propose trades for their own benefits, very much like how today's bilateral power purchase contracts are formed. Thus the fundamental building block of such a system is the notion of trade:
\begin{definition}[Contingent trade]\label{def:ctrade}
	A contingent multilateral trade (referred to as trade in the sequel)\footnote{\rv{See Table~\ref{t:initrade} in Section~\ref{sec:eg} for an example.}} among a group $\Ical^k \subset \Ical$ of participants, is a collection of power injection plans 
	\begin{equation} 
	p^k = (p_{i,s}^k)_{ i\in \Ical^k,\,\, s\in [S]},
	\end{equation}
	that are feasible with respect to participants' local constraints, i.e. $p_i \in \Pcal_i$, and ex post balanced so that
	\begin{equation}
	\sum_{i\in \Ical^k} p^k_{i,s} = 0, \quad s \in [S]. 
	\end{equation} 
\end{definition}
For convenience, we also define $p_{i,s}^k =0$ for $i\not\in \Ical^k$, so that given $p^k$ we can infer $\Ical^k$ via
\begin{equation}
\Ical^k = \{i \in \Ical: \mbox{there exists a $s\in [S]$ s.t. } p^k_{i,s} \neq 0\}. 
\end{equation}

This definition is convenient from the point of view of the SO. In practice,  a trade is a transaction that exchanges power with money. We will touch upon the money side of the trading process in Section~\ref{sec:price} and~\ref{sec:ad}. The power balance condition is natural: the amounts of power supplied and  consumed must be equal in each scenario. 
This  definition also stresses that the commodity for sale is \emph{scenario-contingent power}. That is, 1 MWh in different scenarios of RT are treated as different commodities.  

Some further remarks are in order for Definition~\ref{def:ctrade}.
\begin{remark}[Need for multilateral trades]
	As indicated in Definition~\ref{def:ctrade}, a trade may involve more than two market participants. Although multilateral trades are less common in practice compared with bilateral trades, for our purpose, it is necessary to consider multilateral trades so that the trading process is guaranteed to converge to an efficient dispatch. See \cite{Wu199975} for an example in which bilateral trading fails to converge to the optimal dispatch due to \emph{loop externality} \cite{Hogan92}. When the network does not have cycles, it is possible to show that bilateral trades suffice under certain conditions. See Appendix~\ref{sec:tree}.
\end{remark}

\begin{remark}[SO's sufficient statistics]
	While market participants must keep track of their own power injection plans, from the power system's perspective,  \emph{contingent (network) nodal injection} vectors, calculated from
	\begin{equation}\label{def:qk}
			q_{n,s}^k = \sum_{i \in \Ical_n} p^k_{i,s}, \quad n \in [N], \,\, s\in [S],
	\end{equation}
	 carry all the necessary information for checking the reliability constraints in \eqref{eq:Pcal}. In particular, a trade  among participants  at the same node of network makes zero contribution to the actual nodal injection and thus is not of concern to the SO. 
\end{remark}

Trades motivated by participants interests do not take into account  power system reliability constraints. So it is necessary to have the SO  verify that trades meet the power system  constraints, and in case of violation to \emph{curtail} trades so that compliance is achieved. Throughout the paper, we consider a simple curtailment scheme:
\begin{definition}[Uniform curtailment]\label{def:curt}
	A trade $p^k$ is said to be curtailed if only a portion of the proposed power injection, $\gamma^k p^k$, is accepted by the SO, where $\gamma^k\in [0,1)$  is the curtailment factor  and 
	\begin{equation}
	(\gamma^k p^k)_{i,s} = \gamma^k p^k_{i,s}, \quad i \in \Ical^k,\,\, s\in [S].
	\end{equation}
	For notational convenience, we also define $\gamma^k = 1$ when a trade is accepted without curtailment.\footnote{\rv{Table~\ref{t:cutrade} in Section~\ref{sec:eg} provides an example of a curtailed trade.}}
\end{definition}
\rv{\begin{remark}[Scenario-dependent curtailment]
The uniform curtailment scheme is the simplest curtailment scheme that ensures  \emph{local feasibility} of curtailed trades given that the initial trades satisfy local constraints. That is, given a trade $p^k$ such that $p^k_i \in \Pcal_i$, $i \in \Ical$, the curtailed trade always satisfies $\gamma^k p^k_i \in \Pcal_i$, $i \in \Ical$. It is possible to make the curtailment scenario-dependent, i.e., for each scenario $s\in [S]$, we can pick a different curtailment factor $\gamma^k_s \in [0,1]$. This  curtailment scheme no longer has the local feasibility property if DA participants are involved in the initial trade. In particular, the curtailed trade will not satisfy non-anticipative constraints of DA participants if the curtailment factors for different scenarios are taken to be different values. A hybrid of the uniform curtailment and scenario-dependent curtailment is to use the former when a trade involves DA participants and to  use the latter when it does not. One can verify that all our results hold for this curtailment scheme as well.
\end{remark}
}
During the DA market time window, a sequence of trades will come up for SO's approval. Thus the notion of power system reliability and the calculation of curtailment depend on  trades that are already accepted into the system. We define  a notion of system trading state as follows:
\begin{definition}[Trading state]
	Given a sequence of trades $p^\kappa$ and their curtailment factor $\gamma^\kappa$, $\kappa = 0, \dots, k-1$, the global trading state is the accumulated participants' contingent power injection 
	\begin{equation}
	y^k_{i,s} = \sum_{\kappa = 0}^{k-1}  \gamma^\kappa p^\kappa_{i,s}, \quad i \in \Ical, \,\, s \in [S],
	\end{equation}
	and the network state for the SO is the accumulated network power injection
	\begin{equation}
	x^k_{n,s} = \sum_{\kappa=0}^{k-1} \gamma^\kappa q^\kappa_{n,s} = \sum_{\kappa=0}^{k-1} \sum_{i \in \Ical_n} \gamma^\kappa  p^\kappa_{i,s}, \quad n \in [N],\,\, s \in [S]. 
	\end{equation}
\end{definition}

The network and trading states relate as
\begin{equation}
x_{n,s}^k = \sum_{i \in \Ical_n} y_{i,s}^k. 
\end{equation} 

Given the current system state $x^k$, a characterization for a trade $p^k$ to be feasible for network constraints~\eqref{eq:Pcal} is that its corresponding network injection vector $q^k$ as defined in \eqref{def:qk} satisfies 
\begin{equation}\label{eq:qfeas}
	x^k_s  + q^k_s \in \Pcal, \quad s \in [S]. 
\end{equation}
Define the scenario-contingent feasible set of network injection as
\begin{equation}
\Qcal_s(x_s) = \Pcal - x_s = \{q_s \in \reals^{N}: x_s + q_s \in \Pcal\}, \quad  s\in [S], 
\end{equation}
and $\Qcal(x) = \Qcal_s(x_1) \times \hdots \times \Qcal_S(x_S)$. Then~\eqref{eq:qfeas} is equivalent to $q^k \in \Qcal(x^k)$. 

A potential issue of the trading process, in view of Definition~\ref{def:curt}, is that $\gamma^k$ may have to be $0$ to bring many trades back to feasible. Indeed, if the market participants are proposing trades without any information regarding the current network state $x^k$, then it is likely that many trades overburdening lines which are already congested at $x^k$ will be proposed. To forestall such a possibility, the SO  requires  participants  to only submit trades that are in the \emph{feasible direction} of the network given the current state. 
\begin{definition}[Feasible direction trade]
Given a network state $x^k$, let $\Lcal_s(x^k_s)$ be the set of active (binding) line constraints in scenario $s$,  that is,
\begin{equation}
\Lcal_s(x^k_s) = \{\ell \in [L]: h^\top_\ell x^k_s = \bar f_\ell\}, \quad s \in [S].
\end{equation}
Then a trade $p^k$ is a feasible-direction (FD) trade at $x^k$ if its corresponding network injection $q^k$ as defined in~\eqref{def:qk} satisfies
\begin{equation}\label{eq:fd}
	h^\top_\ell q^k_s \le 0, \quad \ell \in \Lcal_s(x^k_s), \,\, s \in [S]. 
\end{equation}
\end{definition}

If market participants are constrained to propose only FD trades, then it is guaranteed that $\gamma^k >0$ so that every trade updates the network state. At this moment, it is unclear whether such update is favorable in any sense. Formalizing the notion of ``self-interested'' participants, we have the following definition.
\begin{definition}[Worthwhile trade]
	We call a trade $p^k$ an $\epsilon$-worthy trade at trading state $y^k$ if it leads to welfare improvement  no smaller than $\epsilon$, i.e.,
	\begin{equation}\label{eq:worthy}
			\sum_{i \in \Ical^k} U_i(y^k_i + p^k_i) - U_i(y^k_i) \ge \epsilon, 
	\end{equation}
	and an $\epsilon$-unworthy trade if~\eqref{eq:worthy} does not hold. A profitable trade is an $\epsilon$-worthy trade with $\epsilon=0$. 
\end{definition}

Notice that if an $\epsilon$-worthy trade is proposed by some participants and accepted by the SO, then it improves the \emph{social welfare} by at least $\epsilon$ as the power injection plans of participants  not involved in the trade are not changed. 

We can now formalize the coordinated  trading process. 
\begin{enumerate}[Step 1.]
	\item \emph{Initialization.} The SO initializes the system state $x^k$ corresponding to some initial feasible trade $p^k$, $k=0$.
	\item \emph{Announcement.} The SO checks the congestion state of the system at $x^k$, identifies $\Lcal_s(x_s^k)$ for $s\in [S]$ and announces the network loading vectors $h_\ell$, $\ell \in \Lcal_s(x^k_s)$ for each $s\in [S]$. 
	\item  \emph{Trading.} If a profitable trade\footnote{We refer \rv{to} the resulting trading process as an $\epsilon$-trading process when the requirement of profitable trade is replaced by $\epsilon$-worthy trade. \label{ft:etrade} } in the feasible direction $p^k$ is identified, participants arrange it. If no profitable trade is found, go to Step 6.
	\item \emph{Curtailment.} If $p^k$ is not feasible, i.e., the corresponding network injection $q^k$ is such that $q^k \not\in \Qcal(x^k)$, the SO curtails the trade with 
	\begin{equation}
	\gamma^k = \max \left\{\gamma: \gamma q^k \in \Qcal(x^k)\right\} \in (0,1).\end{equation} If $p^k$ is feasible, set $\gamma^k = 1$.
	\item  \emph{Update.} The SO updates the network state as $x^{k+1} \gets x^k + \gamma^k q^k$, $k \gets k+1$. Go to Step 2. 
	\item  \emph{Termination.}
\end{enumerate}

It is evident from the description of the trading process that the SO only has the following responsibilities. (i) SO checks whether the trade newly submitted by participants is feasible with respect to network constraints. If not, it curtails the trade so that the resulting trade is feasible. (ii) In case there are  congested lines, the SO computes and broadcasts the loading vectors to the participants. Note that in our framework, the SO does not carry out any optimization. Instead, market participants  seek to optimize their own profit during the trading process. 

\begin{remark}[Feasibility]
	An important feature of the trading process is that the proposed system state $x^k$ for any $k$ is feasible with respect to the power network constraints. Thus even if the trading process is stopped at any stage before termination, the  trades still result in a safe power flow solution. 
	\end{remark}

\begin{remark}[Pay-as-bid settlement]
The trading process allows a pay-as-bid 	payment settlement approach.  Immediately after submitting the trade, the market participants are informed  whether their trades will be scheduled (or partially scheduled if curtailed); this information can then be used to calculate and settle the payment among these participants. Comparing to the locational pricing used in the standard  market, such a payment settlement process could limit the price risk faced by market participants which are expected to increase when the system integrates more renewables. This is also part of the reason why bilateral long-term contracts are widely used by large utility companies and power producers. 
\end{remark}

\section{Economic Efficiency}\label{sec:eff}
Similar to arguments in  \cite{Wu199975}, one can verify that  trading process described in the previous subsection is well-defined, and whenever a $\epsilon$-worthy trade is identified $(\epsilon >0)$, the social welfare is strictly increased (even if the trade is to be curtailed). Thus when the trading process terminates, that is, when there exists no additional profitable trade that is not yet arranged, one may expect that the resulting power injection plan matches the economic efficiency benchmark defined by stochastic optimization problem~\eqref{opt:swm}.  
\begin{theorem}[Efficiency] \label{thm:main}
Suppose the following assumptions are in force:
\begin{enumerate}[(i)]
\item for any fixed $\epsilon>0$, any $\epsilon$-unworthy trade in the feasible direction will not be arranged and any $\epsilon$-worthy trade will eventually be identified and arranged, and
\item once a worthy profitable trade is identified, the market participants involved are willing to carry it out.
\end{enumerate}
Then the $\epsilon$-trading process  is well-defined and the accumulated global trading state $y^k$ converges in the sense that
\begin{equation}U^\star - \lim_{k \to \infty} U(y^k) \le \epsilon, \end{equation}
 where $U^\star$ is the optimal value of \eqref{opt:swm}.
\end{theorem}

When $\epsilon$ is sufficiently small, Theorem~\ref{thm:main} states that the trading process will converge to a dispatch that is practically optimal for any desired accuracy. 
 The key message of Theorem~\ref{thm:main} is that the extremely simple feedback procedure of the SO based on curtailment and loading vector announcement suffices in providing coordination for the trades so that efficiency is achieved while network reliability is guaranteed in every iteration of the trading process.  
\begin{remark}[Trade formation]
	Like in \cite{Wu199975}, we purposely leave the details of trading group formation open. Theorem~\ref{thm:main} is powerful in that it is agnostic to the actual underlying mechanism dictating which subset of participants meet and propose trade $k$. For instance, a conceptually simple mechanism is that in every iteration $k$, a subset of $\Ical$ is picked  at random such that there is a positive probability for picking every subset\footnote{We can e.g. first sample a random group size $I^k$ from $[|\Ical|]$ and then randomly sample a group from $\Ical$ of size $I^k$. }. If it is possible for this group of participants to identify a profitable trade in the feasible direction, they will propose it, as in Step 3 of the trading process. If not, we can simply continue this process by generating another random group of participants. Since there is a finite number of subsets of $\Ical$, Theorem~\ref{thm:main} guarantees that this process converges to an efficient dispatch with probability one. In practice,  trading group formation  processes depend on a lot of factors that we do not model in this paper. As a result, it could be the case  that each participant $i \in \Ical$ may only have access to a small subset of other participants in the market. An important future research direction is to design \emph{information platforms} that facilitate trade discovery and reduce \emph{search cost}. 
\end{remark}
\begin{remark}[Profit allocation]
	Similarly, we do not specify how profit is allocated  among the participants if a profitable trade is proposed and accepted by the SO. One can verify (or cf. \cite{Wu199975}) that for every profitable trade, there is a profit allocation that makes \emph{all} involving participants better off. 
\end{remark}
\begin{remark}[Merchandising surplus]
	In the standard market,  the total payment collected from loads is larger than that paid to generators when there are line congestions. This merchandising surplus \cite{Wu1996} is paid to transmission owners. In our setting, as the SO  does not collect money from participants and all trades are budget balanced,  separate payment streams might be needed to cover the costs of the transmission owners. Possible ways include charging a fee for using the transmission or requiring participants to acquire transmission rights for making trades across the network. 
\end{remark}

\begin{remark}[Algorithmic interpretation]
	The trading process may be thought of as a \emph{projected line search} algorithm for solving \eqref{opt:swm}, whose iteration $k$ performs update 
	\begin{equation}
	y^{k+1} = y^k + \gamma^k p^k,
	\end{equation}
	where $p^k$ is the search direction and $\gamma^k$ is the step length introduced to project the step into the feasible region. The algorithm is \emph{distributed}, in that the search direction is identified based on  information (and  objective functions) of a subset of participants. The algorithm is special as its search direction $p^k$ is identified from a profitable trade, which is an economic construct, rather than based on gradient or Hessian of the objective function. 
\end{remark}

\rv{\begin{remark}[Subjective probability]
In general, different market participants may have their own subjective assessment of the probabilities of the scenarios. Denote the subjective probabilities  of participant $i \in \mathcal I$ by $\Prob_i(s)$, $s\in [S]$. Then $y^k$ converges to an optimal solution of \eqref{opt:swm} with the ex ante utility function replaced by
\begin{equation}
U_i(p_i) = \E_i [u_{i,s} (p_{i,s})] = \sum_{s\in [S]} \Prob_i(s) u_{i,s}(p_{i,s}).
\end{equation}
In this case, the resulting dispatch is  Pareto optimal but may not maximize the ex ante social welfare as the latter notion is defined upon the unknown true probability distribution of the uncertainty. 
\end{remark}}

\begin{remark}[Distribution system operator]
	The trading process also offers a way to design a lightweight or minimal distribution system operator (minDSO) for coordinating distributed generation (DG), flexible loads and other distribution level resources. With minDSO, the DG owners and demand side flexibility providers do not need to  report cost and benefit data to the minDSO; so long as they can determine profitable trades among themselves, social welfare will improve. To adapt our formulation to the distribution system setting, the linearized DistFlow model\cite{19266, 25627} provides an accurate model of the real power flow on the distribution network. Line capacity and transformer limits can be modeled similarly as transmission line limits. The tree network topology offers potential simplification to the trading process (see Appendix~\ref{sec:tree}). Voltage constraints can be modeled as additional linear constraints \cite{ farivar2013equilibrium}. Distribution topology switching can be accommodated by updating the network constraint set $\Pcal$ according to the current switch states. 
	
\end{remark}

\section{Price Discovery}\label{sec:price}
In standard  markets, the SO solves an economic dispatch problem that determines both the dispatch and the locational marginal prices of power at all buses. When uncertainty is considered,  the computationally demanding stochastic economic dispatch must be solved by the SO. 

The trading process, on the contrary, does not require the SO to solve any optimization problem. Theorem~\ref{thm:main} suggests that an efficient dispatch is achieved in the limit; here we show that the optimal locational marginal prices also \emph{emerge} when the trading process converges\footnote{An alternative treatment, involving setting up trade-based prices and characterizing the convergence of the price process, is also possible (cf. \cite{Chao1996,RePEc:eee:mateco:v:63:y:2016:i:c:p:84-92}). However, this requires a detailed specification of  payments associated with each trade which we avoid in this paper.}. The idea is simple. Suppose that in the last few minutes of the DA trading window, when the trading process has already converged, a new load comes into the system and demands a $\epsilon\to 0$ unit of power at bus $n$ for scenario $s$.  Producers who can still generate additional power could  each quote a price  based on their \emph{marginal cost} evaluated at the current trading state. We thus discover the locational marginal price at bus $n$ for scenario $s$ by finding the minimum price announced by  those generators that can indeed send power to bus $n$ given the congestion state in the scenario. 

To formalize this idea, denote the optimal dual variable associated with constraint~\eqref{opt:swm:pb} by $\lambda_{n,s}^\star$, $n \in [N]$, $s\in [S]$. Furthermore, as constraint~\eqref{opt:swm:local} for RT participants is a box constraint, denote the optimal dual variables associated with the lower and upper bounds by $\underline{\eta}_{i,s}^\star$ and $\overline{\eta}_{i,s}^\star$, respectively. Notice that the trading process has a balanced budget by construction, as the system operator is not involved in any financial aspect of the system. Therefore, problem~\eqref{opt:swm} is equivalent to the stochastic economic dispatch problem. 
\begin{lemma}[Price discovery]\label{lemma:price}
	For each bus $n \in [N]$ and $s\in [S]$, if there exists a participant $i\in \IRT_n$ whose utility function is differentiable and whose optimal contingent power injection $p^\star_{i,s}$ is in the interior of her local feasible set, i.e., $p^\star_{i,s} \in \mathring{\Pcal_{i,s}}$, then we have\footnote{The sign convention is  that $p_{i,s}>0$ represents power injection (supply) into the network. Thus if $\tilde u_{i,s}$ is a utility function that in the usual sense is increasing with  demand,  $\lambda^\star_{n,s} $ as computed below will be nonnegative. }
	\begin{equation}\label{eq:lmp:int}
		\lambda^\star_{n,s} = - \rev{\Prob}(s) \frac{\partial \tilde u_{i,s}(p^\star_{i,s})}{\partial p_{i,s}}. 
	\end{equation}  
	In general, suppose the utility function of some participant $i \in \IRT_n$ is differentiable, then 
	\begin{equation}\label{eq:lmp:gen}
			\lambda^\star_{n,s} = - \rev{\Prob}(s) \frac{\partial \tilde u_{i,s}(p^\star_{i,s})}{\partial p_{i,s}} + (\overline{\eta}_{i,s}^\star - \underline{\eta}_{i,s}^\star).
	\end{equation}
\end{lemma}

While the price calculation based on~\eqref{eq:lmp:int} is intuitive and only requires local information, that based on~\eqref{eq:lmp:gen} may require solving the dual program of~\eqref{opt:swm} to identify the values of the optimal dual variables $\overline{\eta}_{i,s}^\star$ and $ \underline{\eta}_{i,s}^\star$. Fortunately, solving for the dual program is greatly simplified when the optimal primal solution $p^\star$ is known (cf. \cite{Boyd:2004:CO:993483}). 


\section{Arrow-Debreu Equilibrium}\label{sec:ad}
Section~\ref{sec:eff} established that the trading process converges to a stationary contingent power injection plan $p^\star$. Section~\ref{sec:price} then showed that there is a well-defined notion of price $\lambda^\star$ that emerges alongside with the stationary injection plan. Here we connect the pair $(p^\star, \lambda^\star)$ to the suitable economic concept of general equilibrium under uncertainty. Taken together, this will formally establish that the contingent trading process converges to a \emph{(contingent plan, price)} equilibrium, which respects the power system reliability constraints and achieves economic efficiency. 

To start, we need to define an electricity market economy, similar to that done in \cite{7447005} (also see \cite{Wang20114933}). The commodities of the economy is \emph{contingent power} at each node $n \in [N]$ and in each scenario $s\in [S]$. Buying (selling) a unit contingent power $(n,s)$ in DA leads to the right to consume (responsibility to generate) a unit of power at node $n$ if scenario $s$ occurs in RT. 

The market participants are those in $\Ical$ as defined in Section~\ref{sec:part} and a \emph{traditional system operator}\footnote{This notion of SO is consistent with that in the literature and different from the SO described in Section~\ref{sec:trading}.} who may convert power at one node into that at another node using the network. 

For each participant $i \in \Ical$, given prices for contingent power $\lambda$, the following optimization is solved to determine the participant's contingent power injection plan
\begin{subequations}\label{opt:i}
	\begin{align}
		\mbox{maximize} \quad & \sum_{s \in [S]} \lambda_{n,s} p_{i,s} + \rev{\Prob}(s)\tilde u_{i,s}(p_{i,s})\label{opt:i:obj}\\
		\mbox{subject to} \quad & p_{i} \in \Pcal_{i}, 
	\end{align}
\end{subequations}
where $\lambda_{n,s}$ is the price at node $n$ faced by $i \in \Ical_n$. 
Here the objective function is the same as $U_i(p_i)$ as defined in~\eqref{eq:vnmU} with linear payment scheme
\begin{equation}
m_{i,s}(p_{i,s}) = \tilde \lambda_{n,s} p_{i,s},
\end{equation}
where \rev{$\tilde \lambda_{n,s}  =\lambda_{n,s}/\Prob(s)$}. Notice that the first term in the summation in \eqref{opt:i:obj} is the monetary payment that clears in DA; the second term is the expected utility derived from the power injection in RT. 

The SO is modeled as a firm that uses technology (i.e. power network) to convert one type of commodity (i.e. contingent power on one node) to other types of commodities (i.e. contingent power on other nodes), in order to maximize its profit. Formally, the SO solves the following optimization to determine the contingent network power injection $x$ given prices for contingent power $\lambda$:
\begin{subequations}\label{opt:so}
\begin{align}
	\mbox{maximize} \quad & \sum_{s\in [S]}\sum_{n\in [N]} -\lambda_{n,s} x_{n,s} \label{opt:so:obj}\\
	\mbox{subject to}\quad & x_s \in \Pcal, \quad  s\in [S],
\end{align}
\end{subequations}
where the entire profit of SO in \eqref{opt:so:obj} is cleared in DA. 

The suitable notion of competitive equilibrium of such a market for contingent claims is that of \emph{Arrow-Debreu}, which generalizes the Walrasian concept of general equilibrium to settings with uncertainty \cite{mas1995microeconomic}. Here we state the definition of Arrow-Debreu equilibrium for the electricity market economy:
\begin{definition}[Arrow-Debreu equilibrium]
	A collection of contingent power injection plans $(p^\star,x^\star)$, with $p^\star\in \reals^{|\Ical|\times S}$ and $x^\star\in \reals^{N\times S}$ and a system of prices for contingent power $\lambda^\star \in \reals^{N\times S}$ constitute an Arrow-Debreu equilibrium if:
	\begin{enumerate}[(i)]
		\item For every $i \in \Ical$, $p_i^\star$ solves \eqref{opt:i} given prices $\lambda^\star$.
		\item For the SO, $x^\star$ solves \eqref{opt:so} given prices $\lambda^\star$.
		\item The  market for each contingent power commodity clears:
		\begin{equation}\label{eq:req}
					 x_{n,s}^\star = \sum_{i \in \Ical_n} p_{i,s}^\star, \quad n\in [N],\,\, s \in [S].
		\end{equation}

	\end{enumerate}
\end{definition}

By the first fundamental theorem of welfare economics and in a quasi-linear environment, we expect that a dispatch-price tuple $(p^\star,x^\star,\lambda^\star)$ at an  Arrow-Debreu equilibrium achieves economic efficiency defined by \eqref{opt:swm} (cf. \cite{mas1995microeconomic, NET-011}). Our previous result suggests that the dispatch at the limit of the trading process together with the emerged prices matches the solution of \eqref{opt:swm}. Our next result establishes that the dispatch-price tuple obtained from the trading process indeed constitutes an Arrow-Debreu equilibrium.
\begin{lemma}\label{lemma:eq}
	Suppose that an Arrow-Debreu equilibrium exists and that the utility functions are differentiable. Then the contingent power injection plan $p^\star$ obtained from the trading process, the corresponding network injection plan $x^\star$ calculated from \eqref{opt:swm:pb}, and the prices computed from~\eqref{eq:lmp:int} or \eqref{eq:lmp:gen} constitute an Arrow-Debreu equilibrium. 
\end{lemma}
\begin{remark}[T\^{a}tonnement process]\label{rk:tp}
	In light of Lemma~\ref{lemma:eq}, the trading process can be thought of as a way to drive an out-of-equilibrium market into its equilibrium. Such processes, characterizing the dynamic laws of out-of-equilibrium movement of the market state, is in general referred to as a t\^{a}tonnement process; see \cite{hahn1982stability}.
\end{remark}

\section{Examples} \label{sec:eg}
We provide an illustrative example for the trading process in this section.

\begin{figure}[htbp]
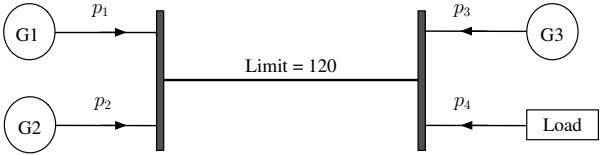

\centering
\scalebox{.45}{\scalefont{2} 
\begin{pgfpicture}
\pgfpathmoveto{\pgfqpoint{0cm}{0cm}}
\pgfpathlineto{\pgfqpoint{18.344cm}{0cm}}
\pgfpathlineto{\pgfqpoint{18.344cm}{5.433cm}}
\pgfpathlineto{\pgfqpoint{0cm}{5.433cm}}
\pgfpathclose
\pgfusepath{clip}
\begin{pgfscope}
\end{pgfscope}
\begin{pgfscope}
\begin{pgfscope}
\pgfpathmoveto{\pgfqpoint{0cm}{0cm}}
\pgfpathlineto{\pgfqpoint{18.344cm}{0cm}}
\pgfpathlineto{\pgfqpoint{18.344cm}{5.433cm}}
\pgfpathlineto{\pgfqpoint{0cm}{5.433cm}}
\pgfpathclose
\pgfusepath{clip}
\definecolor{eps2pgf_color}{rgb}{1,1,1}\pgfsetstrokecolor{eps2pgf_color}\pgfsetfillcolor{eps2pgf_color}
\pgfpathmoveto{\pgfqpoint{-0.723cm}{6.262cm}}
\pgfpathlineto{\pgfqpoint{39.917cm}{6.262cm}}
\pgfpathlineto{\pgfqpoint{39.917cm}{-19.597cm}}
\pgfpathlineto{\pgfqpoint{-0.723cm}{-19.597cm}}
\pgfpathclose
\pgfusepath{fill}
\pgfsetdash{}{0cm}
\pgfsetlinewidth{0.706mm}
\pgfsetroundcap
\pgfsetroundjoin
\definecolor{eps2pgf_color}{rgb}{0,0,0}\pgfsetstrokecolor{eps2pgf_color}\pgfsetfillcolor{eps2pgf_color}
\pgfpathmoveto{\pgfqpoint{5.115cm}{2.54cm}}
\pgfpathlineto{\pgfqpoint{12.629cm}{2.54cm}}
\pgfusepath{stroke}
\pgfsetdash{}{0cm}
\pgfsetlinewidth{0.353mm}
\pgfpathmoveto{\pgfqpoint{1.683cm}{4.555cm}}
\pgfpathcurveto{\pgfqpoint{1.979cm}{4.231cm}}{\pgfqpoint{1.979cm}{3.706cm}}{\pgfqpoint{1.683cm}{3.383cm}}
\pgfpathcurveto{\pgfqpoint{1.387cm}{3.059cm}}{\pgfqpoint{0.906cm}{3.059cm}}{\pgfqpoint{0.61cm}{3.383cm}}
\pgfpathcurveto{\pgfqpoint{0.314cm}{3.706cm}}{\pgfqpoint{0.314cm}{4.231cm}}{\pgfqpoint{0.61cm}{4.555cm}}
\pgfpathcurveto{\pgfqpoint{0.906cm}{4.879cm}}{\pgfqpoint{1.387cm}{4.879cm}}{\pgfqpoint{1.683cm}{4.555cm}}
\pgfusepath{stroke}
\definecolor{eps2pgf_color}{rgb}{0,0,0}\pgfsetstrokecolor{eps2pgf_color}\pgfsetfillcolor{eps2pgf_color}
\pgftext[x=1.086cm,y=3.889cm,rotate=0]{\fontsize{16}{11.76}\selectfont{G1}}
\definecolor{eps2pgf_color}{rgb}{0.317884,0.317825,0.317864}\pgfsetstrokecolor{eps2pgf_color}\pgfsetfillcolor{eps2pgf_color}
\pgfpathmoveto{\pgfqpoint{4.904cm}{4.586cm}}
\pgfpathlineto{\pgfqpoint{5.115cm}{4.586cm}}
\pgfpathlineto{\pgfqpoint{5.115cm}{0.459cm}}
\pgfpathlineto{\pgfqpoint{4.904cm}{0.459cm}}
\pgfpathclose
\pgfusepath{fill}
\pgfsetdash{}{0cm}
\definecolor{eps2pgf_color}{rgb}{0,0,0}\pgfsetstrokecolor{eps2pgf_color}\pgfsetfillcolor{eps2pgf_color}
\pgfpathmoveto{\pgfqpoint{4.904cm}{4.586cm}}
\pgfpathlineto{\pgfqpoint{5.115cm}{4.586cm}}
\pgfpathlineto{\pgfqpoint{5.115cm}{0.459cm}}
\pgfpathlineto{\pgfqpoint{4.904cm}{0.459cm}}
\pgfpathclose
\pgfusepath{stroke}
\definecolor{eps2pgf_color}{rgb}{0.317884,0.317825,0.317864}\pgfsetstrokecolor{eps2pgf_color}\pgfsetfillcolor{eps2pgf_color}
\pgfpathmoveto{\pgfqpoint{12.629cm}{4.586cm}}
\pgfpathlineto{\pgfqpoint{12.841cm}{4.586cm}}
\pgfpathlineto{\pgfqpoint{12.841cm}{0.459cm}}
\pgfpathlineto{\pgfqpoint{12.629cm}{0.459cm}}
\pgfpathclose
\pgfusepath{fill}
\pgfsetdash{}{0cm}
\definecolor{eps2pgf_color}{rgb}{0,0,0}\pgfsetstrokecolor{eps2pgf_color}\pgfsetfillcolor{eps2pgf_color}
\pgfpathmoveto{\pgfqpoint{12.629cm}{4.586cm}}
\pgfpathlineto{\pgfqpoint{12.841cm}{4.586cm}}
\pgfpathlineto{\pgfqpoint{12.841cm}{0.459cm}}
\pgfpathlineto{\pgfqpoint{12.629cm}{0.459cm}}
\pgfpathclose
\pgfusepath{stroke}
\pgfsetdash{}{0cm}
\pgfsetlinewidth{0.423mm}
\pgfpathmoveto{\pgfqpoint{1.922cm}{3.945cm}}
\pgfpathlineto{\pgfqpoint{2.176cm}{3.945cm}}
\pgfpathlineto{\pgfqpoint{4.904cm}{3.945cm}}
\pgfpathlineto{\pgfqpoint{4.904cm}{3.854cm}}
\pgfusepath{stroke}
\pgfsetdash{}{0cm}
\pgfpathmoveto{\pgfqpoint{1.923cm}{1.199cm}}
\pgfpathlineto{\pgfqpoint{4.904cm}{1.192cm}}
\pgfusepath{stroke}
\pgfsetdash{}{0cm}
\pgfpathmoveto{\pgfqpoint{12.841cm}{3.986cm}}
\pgfpathlineto{\pgfqpoint{15.822cm}{3.979cm}}
\pgfusepath{stroke}
\pgfsetdash{}{0cm}
\pgfpathmoveto{\pgfqpoint{12.841cm}{1.199cm}}
\pgfpathlineto{\pgfqpoint{15.822cm}{1.192cm}}
\pgfusepath{stroke}
\pgfsetdash{}{0cm}
\pgfsetlinewidth{0.353mm}
\pgfpathmoveto{\pgfqpoint{15.84cm}{1.658cm}}
\pgfpathlineto{\pgfqpoint{17.956cm}{1.658cm}}
\pgfpathlineto{\pgfqpoint{17.956cm}{0.776cm}}
\pgfpathlineto{\pgfqpoint{15.84cm}{0.776cm}}
\pgfpathclose
\pgfusepath{stroke}
\definecolor{eps2pgf_color}{rgb}{0,0,0}\pgfsetstrokecolor{eps2pgf_color}\pgfsetfillcolor{eps2pgf_color}
\pgftext[x=16.899cm,y=1.217cm,rotate=0]{\fontsize{16}{11.76}\selectfont{Load}}
\pgftext[x=8.881cm,y=2.963cm,rotate=0]{\fontsize{16}{11.76}\selectfont{Limit = 120}}
\pgfsetdash{}{0cm}
\definecolor{eps2pgf_color}{rgb}{0,0,0}\pgfsetstrokecolor{eps2pgf_color}\pgfsetfillcolor{eps2pgf_color}
\pgfpathmoveto{\pgfqpoint{1.7cm}{1.803cm}}
\pgfpathcurveto{\pgfqpoint{1.997cm}{1.48cm}}{\pgfqpoint{1.997cm}{0.955cm}}{\pgfqpoint{1.7cm}{0.631cm}}
\pgfpathcurveto{\pgfqpoint{1.404cm}{0.307cm}}{\pgfqpoint{0.924cm}{0.307cm}}{\pgfqpoint{0.628cm}{0.631cm}}
\pgfpathcurveto{\pgfqpoint{0.332cm}{0.955cm}}{\pgfqpoint{0.332cm}{1.48cm}}{\pgfqpoint{0.628cm}{1.803cm}}
\pgfpathcurveto{\pgfqpoint{0.924cm}{2.127cm}}{\pgfqpoint{1.404cm}{2.127cm}}{\pgfqpoint{1.7cm}{1.803cm}}
\pgfusepath{stroke}
\definecolor{eps2pgf_color}{rgb}{0,0,0}\pgfsetstrokecolor{eps2pgf_color}\pgfsetfillcolor{eps2pgf_color}
\pgftext[x=1.16cm,y=1.14cm,rotate=0]{\fontsize{16}{11.76}\selectfont{G2}}
\pgfsetdash{}{0cm}
\definecolor{eps2pgf_color}{rgb}{0,0,0}\pgfsetstrokecolor{eps2pgf_color}\pgfsetfillcolor{eps2pgf_color}
\pgfpathmoveto{\pgfqpoint{17.135cm}{4.555cm}}
\pgfpathcurveto{\pgfqpoint{17.431cm}{4.231cm}}{\pgfqpoint{17.431cm}{3.706cm}}{\pgfqpoint{17.135cm}{3.383cm}}
\pgfpathcurveto{\pgfqpoint{16.838cm}{3.059cm}}{\pgfqpoint{16.358cm}{3.059cm}}{\pgfqpoint{16.062cm}{3.383cm}}
\pgfpathcurveto{\pgfqpoint{15.766cm}{3.706cm}}{\pgfqpoint{15.766cm}{4.231cm}}{\pgfqpoint{16.062cm}{4.555cm}}
\pgfpathcurveto{\pgfqpoint{16.358cm}{4.879cm}}{\pgfqpoint{16.838cm}{4.879cm}}{\pgfqpoint{17.135cm}{4.555cm}}
\pgfusepath{stroke}
\definecolor{eps2pgf_color}{rgb}{0,0,0}\pgfsetstrokecolor{eps2pgf_color}\pgfsetfillcolor{eps2pgf_color}
\pgftext[x=16.593cm,y=3.891cm,rotate=0]{\fontsize{16}{11.76}\selectfont{G3}}
\pgfsetdash{}{0cm}
\definecolor{eps2pgf_color}{rgb}{0,0,0}\pgfsetstrokecolor{eps2pgf_color}\pgfsetfillcolor{eps2pgf_color}
\pgfpathmoveto{\pgfqpoint{3.228cm}{3.951cm}}
\pgfpathlineto{\pgfqpoint{3.69cm}{3.951cm}}
\pgfusepath{stroke}
\pgfpathmoveto{\pgfqpoint{3.972cm}{3.951cm}}
\pgfpathlineto{\pgfqpoint{3.69cm}{4.057cm}}
\pgfpathlineto{\pgfqpoint{3.69cm}{3.845cm}}
\pgfpathclose
\pgfusepath{fill}
\pgfsetdash{}{0cm}
\pgfsetbuttcap
\pgfsetmiterjoin
\pgfpathmoveto{\pgfqpoint{3.972cm}{3.951cm}}
\pgfpathlineto{\pgfqpoint{3.69cm}{4.057cm}}
\pgfpathlineto{\pgfqpoint{3.69cm}{3.845cm}}
\pgfpathclose
\pgfusepath{stroke}
\definecolor{eps2pgf_color}{rgb}{0,0,0}\pgfsetstrokecolor{eps2pgf_color}\pgfsetfillcolor{eps2pgf_color}
\pgftext[x=3.255cm,y=4.582cm,rotate=0]{\fontsize{16}{11.76}\selectfont{$p_1$}}
\pgfsetdash{}{0cm}
\pgfsetroundcap
\pgfsetroundjoin
\definecolor{eps2pgf_color}{rgb}{0,0,0}\pgfsetstrokecolor{eps2pgf_color}\pgfsetfillcolor{eps2pgf_color}
\pgfpathmoveto{\pgfqpoint{3.237cm}{1.199cm}}
\pgfpathlineto{\pgfqpoint{3.699cm}{1.199cm}}
\pgfusepath{stroke}
\pgfpathmoveto{\pgfqpoint{3.981cm}{1.199cm}}
\pgfpathlineto{\pgfqpoint{3.699cm}{1.305cm}}
\pgfpathlineto{\pgfqpoint{3.699cm}{1.094cm}}
\pgfpathclose
\pgfusepath{fill}
\pgfsetdash{}{0cm}
\pgfsetbuttcap
\pgfsetmiterjoin
\pgfpathmoveto{\pgfqpoint{3.981cm}{1.199cm}}
\pgfpathlineto{\pgfqpoint{3.699cm}{1.305cm}}
\pgfpathlineto{\pgfqpoint{3.699cm}{1.094cm}}
\pgfpathclose
\pgfusepath{stroke}
\definecolor{eps2pgf_color}{rgb}{0,0,0}\pgfsetstrokecolor{eps2pgf_color}\pgfsetfillcolor{eps2pgf_color}
\pgftext[x=3.321cm,y=1.832cm,rotate=0]{\fontsize{16}{11.76}\selectfont{$p_2$}}
\pgfsetdash{}{0cm}
\pgfsetroundcap
\pgfsetroundjoin
\definecolor{eps2pgf_color}{rgb}{0,0,0}\pgfsetstrokecolor{eps2pgf_color}\pgfsetfillcolor{eps2pgf_color}
\pgfpathmoveto{\pgfqpoint{14.205cm}{3.986cm}}
\pgfpathlineto{\pgfqpoint{14.667cm}{3.986cm}}
\pgfpathlineto{\pgfqpoint{14.667cm}{3.986cm}}
\pgfusepath{stroke}
\pgfpathmoveto{\pgfqpoint{13.923cm}{3.986cm}}
\pgfpathlineto{\pgfqpoint{14.205cm}{3.881cm}}
\pgfpathlineto{\pgfqpoint{14.205cm}{4.092cm}}
\pgfpathclose
\pgfusepath{fill}
\pgfsetdash{}{0cm}
\pgfsetbuttcap
\pgfsetmiterjoin
\pgfpathmoveto{\pgfqpoint{13.923cm}{3.986cm}}
\pgfpathlineto{\pgfqpoint{14.205cm}{3.881cm}}
\pgfpathlineto{\pgfqpoint{14.205cm}{4.092cm}}
\pgfpathclose
\pgfusepath{stroke}
\definecolor{eps2pgf_color}{rgb}{0,0,0}\pgfsetstrokecolor{eps2pgf_color}\pgfsetfillcolor{eps2pgf_color}
\pgftext[x=13.939cm,y=4.601cm,rotate=0]{\fontsize{16}{11.76}\selectfont{$p_3$}}
\pgfsetdash{}{0cm}
\pgfsetroundcap
\pgfsetroundjoin
\definecolor{eps2pgf_color}{rgb}{0,0,0}\pgfsetstrokecolor{eps2pgf_color}\pgfsetfillcolor{eps2pgf_color}
\pgfpathmoveto{\pgfqpoint{14.214cm}{1.199cm}}
\pgfpathlineto{\pgfqpoint{14.676cm}{1.199cm}}
\pgfpathlineto{\pgfqpoint{14.676cm}{1.199cm}}
\pgfusepath{stroke}
\pgfpathmoveto{\pgfqpoint{13.932cm}{1.199cm}}
\pgfpathlineto{\pgfqpoint{14.214cm}{1.094cm}}
\pgfpathlineto{\pgfqpoint{14.214cm}{1.305cm}}
\pgfpathclose
\pgfusepath{fill}
\pgfsetdash{}{0cm}
\pgfsetbuttcap
\pgfsetmiterjoin
\pgfpathmoveto{\pgfqpoint{13.932cm}{1.199cm}}
\pgfpathlineto{\pgfqpoint{14.214cm}{1.094cm}}
\pgfpathlineto{\pgfqpoint{14.214cm}{1.305cm}}
\pgfpathclose
\pgfusepath{stroke}
\definecolor{eps2pgf_color}{rgb}{0,0,0}\pgfsetstrokecolor{eps2pgf_color}\pgfsetfillcolor{eps2pgf_color}
\pgftext[x=13.952cm,y=1.832cm,rotate=0]{\fontsize{16}{11.76}\selectfont{$p_4$}}
\end{pgfscope}
\end{pgfscope}
\end{pgfpicture}}
\caption{Network diagram for the two-bus example.}
\label{fig:twobus}
\end{figure}

Consider a two bus network depicted in Figure~\ref{fig:twobus}. There are 3 generators and 1 load connected to the system. We list the relevant data for the participants as follows
\begin{itemize}
\item G1 is a coal power plant that can generate up to $200$ MW at a constant marginal cost $50$ \$/MW. It can only be scheduled in the DA stage due to its lead time.
\item G2 is a wind farm that generates $100$ MW in the first scenario (windy scenario) and $50$ MW in the second scenario (breezy scenario) at no operational cost. Suppose there are only these two scenarios for the system and one of them is realized at the delivery time. The underlying probabilities for these two scenarios are $0.6$ and $0.4$, respectively. 
\item G3 is a gas power plant that can ramp up rapidly at real time with $100$ MW capacity and constant marginal cost of $80$ \$/MW.
\item Load represents an inelastic power consumption of $150$ MW.
\end{itemize}

For the purpose of illustrating the interaction between market participants and the PSO, in this example we assume that in each iteration all four participants meet and propose a trade with no knowledge of the network constraint. 
In DA, as an example, the participants could solve the following optimization problem to identify the cost minimization trade
\begin{align*}
\mbox{minimize} \quad & 50 p_{1} + \E \left[80 p_{3,s}\right] =50 p_1 + 48 p_{3,1} + 32 p_{3,2}\\
\mbox{subject to} \quad &  p_{1} + p_{2,s} + p_{3,s} +p_{4} = 0, \quad s =1,2,\\
                        &  p_4 = -150,\\
                        & 0\le p_{1} \le 200,\\
                        & 0 \le p_{3,s} \le 100,\quad s = 1,2,\\
                        & 0 \le p_{2,1} \le 100, \quad 0 \le p_{2,2} \le 50,
\end{align*}
where the optimization variables are the day-ahead scheduled coal power generation $p_1$, the real-time gas power generation $p_{3,s}$ corresponding to the two scenarios, and wind power generation corresponding to the two scenarios $p_{2,s}$ (which is controllable up to curtailment). 
Upon solving this linear program, the participants propose its solution as their initial trade to the SO, which is shown in Table~\ref{t:initrade}.
\begin{table}[hbt]
 \caption{Power injection (unit: MW) of the initial trade proposed by the particpants.} \label{t:initrade}
\centering
  \begin{tabular}{l|cccc}\toprule
    Scenario & G1 & G2 & G3 & Load\\
    \midrule
    Windy & 50& 100 & 0  & -150\\
    Breezy& 50& 50 & 50 & -150\\\bottomrule
  \end{tabular}
\end{table}

This trade is not feasible with respect to the line limit in the windy scenario. As such, the SO curtails the trade to the one shown in Table~\ref{t:cutrade} with $\gamma = 0.8$. The SO also announces the loading vector such that the constraints~\eqref{eq:fd} can be expressed as $\Delta p_{1}+ \Delta p_{2,s}  -\Delta p_{3,s} - \Delta p_{4} \le 0 $, where $\Delta p$'s are the corresponding changes in the power injections, and $s = 1,2$ as the line limit constraint is binding for both scenarios. The participants then solve the following program to identify a profitable trade in the feasible direction:
\begin{table}[hbt]
\caption{Power injection (unit: MW) of the curtailed trade.} \label{t:cutrade}
\centering
  \begin{tabular}{l|cccc}\toprule
    Scenario & G1 & G2 & G3 & Load\\
    \midrule
    Windy & 40 & 80 & 0  & -120\\
    Breezy& 40& 40 & 40 & -120\\\bottomrule
  \end{tabular}
  \end{table}
\begin{align*}
\mbox{minimize} \,\, &  50 (\tilde p_1\!+\!\Delta p_1) \!+\! 48 (\tilde p_{3,1} \!+\! \Delta p_{3,1}) \!+\! 32 (\tilde p_{3,2} \!+\!\Delta p_{3,2})\\
\mbox{subject to} \,\, &  \Delta p_{1} + \Delta p_{2,s} + \Delta p_{3,s} + \Delta p_4 = 0,\quad s = 1,2, \\
                        & \Delta p_{1} + \Delta p_{2,s}  -\Delta p_{3,s} - \Delta p_{4} \le 0, \quad s = 1,2,\\
                        &  \tilde p_4 + \Delta p_4 = -150,\\
                        & 0\le \tilde p_{1} + \Delta p_1 \le 200,\\
                        & 0 \le \tilde p_{3,s} + \Delta p_{3,s} \le 100,\quad s = 1,2,\\
                        & 0 \le \tilde p_{2,1} \!+\! \Delta p_{2,1} \le 100, \quad 0 \le \tilde p_{2,2} \!+\! \Delta p_{2,2} \le 50,
\end{align*}
where the $\tilde p$'s are the curtailed trade given in Table~\ref{t:cutrade}. The resulting \emph{accumulated} trade $\gamma p+ \Delta p$ is shown in Table~\ref{t:result}. The trading process would terminate now as there is no profitable feasible direction trade can be further identified. One can easily verify that the accumulated trade coincides with the solution to~\eqref{opt:swm} for this example, and therefore the trading process indeed achieves efficiency.
\begin{table}[hbt]
\caption{Power injection (unit: MW) of the accumulated trade $\gamma p + \Delta p$.}
\centering
  \begin{tabular}{l|cccc}\toprule
    Scenario & G1 & G2 & G3 & Load\\
    \midrule
    Windy & 20 & 100 & 30  & -150\\
    Breezy& 20&  50& 80 & -150\\\bottomrule
  \end{tabular}
  \label{t:result}
\end{table}

\section{Concluding remarks and open questions} \label{sec:conclusion}
Contingent coordinated  trading is proposed as a market framework for power system resource allocation under uncertainty. Within the framework, the economic efficiency is achieved via coordinated  trades proposed by any groups of market participants for their own benefit. The trading process also discovers the optimal contingent locational marginal prices, and supports an Arrow-Debreu equilibrium of the market. 
Allowing the trades to be contingent on properly defined system scenarios greatly enhances the flexibility of the trades and could result in an improvement of social welfare compared to standard deterministic dispatch based market clearing. The role of the SO is minimal in our framework as  the SO only monitors the trades, curtails them if necessary, and does not collect any cost data or directly intervene in economic decisions. As such, all suppliers and consumers have open access to the power network, which promotes competition and expedites the processes of new generation and consumer-side technology adoption.

\rv{We envision that the proposed framework could help address  many challenges in designing  new DSOs for distribution systems with deep distributed energy resource penetration. Given the novelty of the proposed framework, it is natural that this paper leaves a variety of fundamental questions open.
\begin{itemize}
	\item \emph{Uncertainty model}. In practice, it is unlikely that we can obtain an exact characterization of all possible scenarios for the entire system in DA. Thus extending our ideal uncertainty model by incorporating information updates could make the trading framework more realistic. Under such settings, it may become advantageous to allow RT re-trading as the realized RT scenario may not be exactly one of the pre-scribed scenarios in DA. Additionally, even if it is possible to characterize the set of all possible scenarios, the total number of scenarios may be very large due to the fact that many scenarios are local (see Appendix~\ref{app:localscenario} for a model where all scenarios are local). Thus in practice, suitable factorization (decomposing the scenario tree into system wide scenarios and local scenarios) or scenario reduction is necessary for successful market design based on the proposed trading framework.
	\item \emph{Trade implementation}. As all  trades happen before real time, in real time, the participants need to supply and consume according to the scheduled trades. To ensure this indeed happens, advance metering infrastructure (AMI) systems and suitable financial incentive (or penalty) scheme have to be in place. Thus an open question is how to design such financial schemes that encourage consistent participant behaviors while limits potential gaming activities. 
	\item  \emph{Trade formation}. For distribution system applications, requiring participants to meet and trade seems overwhelming. A more likely setting is  to rely on one or many third-party marketplaces to identify  profitable trades on behalf of (subsets of) the participants. Our analysis also  applies to such settings thanks to our general assumption on the trade formation process.  In this context, our results are better understood as a form of \emph{separation principle}, which ensures lossless separation of network reliability from market efficiency considerations  with our trading framework. Under such separation, third party marketplaces can fill the role of trade identification and formation without any explicit knowledge of the power network, as long as it follows the rules set by the SO. 
	 Designing and implementing such third-party marketplaces to unlock potentials from distributed energy resources and flexible loads thus is an important future direction to explore. 
\end{itemize}}

\section*{Acknowledgement}
The authors are grateful to Shmuel Oren, Kameshwar Poolla, Eilyan Bitar and Subhonmesh Bose for very helpful inputs and discussions.  

\bibliographystyle{IEEEtran}
\bibliography{OA1}

\appendices
\section{Proof of Theorem \ref{thm:main}}
We prove Theorem~\ref{thm:main} by first establishing a more general version of the result for optimization of the form 
 \begin{subequations}\label{opt:abs}
	\begin{align}
		\mbox{maximize} \quad & U(p)\\
		\mbox{subject to} \quad & Ep = 0,\\
		& Dp \le d, \\
		& p_i \in \Pcal_i, \quad i \in [m],
	\end{align}
\end{subequations}
and then show that Theorem~\ref{thm:main} is a special case. 
Here the decision variables  are $p_i \in \reals^n$, $i \in [m]$ so that $p \in \reals^{mn}$, constraint matrices   $E \in \reals^{r_1 \times mn}$ and $D\in \reals^{r_2 \times mn}$ are such that $\{p: Ep=0, Dp \le d\} $ is a closed compact  subset of $\reals^{mn}$, local feasible sets $\mathcal P_i$, $i \in [m]$ are  closed, compact and convex subsets of $\reals^n$, where $n, m, r_1$ and $r_2$ are positive integers. 

For each $\epsilon>0$, we consider an iterative algorithm for solving \eqref{opt:abs} of the following form: given $y^0 \in \reals^{mn}$ feasible for \eqref{opt:abs}, the iterations are generated by a point-to-set mapping $A_\epsilon: \reals^{mn}\mapsto 2^{\reals^{mn}}$. That is, 
\begin{equation}\label{eq:alg}
	y^{k+1} \in A_\epsilon(y^k), \quad k \in \Zbb_+. 
\end{equation}
We consider a specific algorithmic  mapping that is a composition of two mappings: $A_\epsilon = CF_\epsilon$. Here the point-to-set mapping $F_\epsilon: \reals^{mn} \mapsto \reals^{mn} \times 2^{\reals^{mn}}$ is defined such that $(y',p) \in F_\epsilon(y)$ if $y' = y$ and 
\begin{equation}
p \in G_\epsilon(y) = \Set{ p \in \reals^{mn}: \begin{split}
\Ical^k \subset [m], p_i =0, \,\, i \not \in \Ical^k,\\
 	y_i + p_i \in \Pcal_i, \,\, i \in \Ical^k,\\
  Ep =0, \,\, D(y)p \le 0,\\
  U(y+ p) - U(y) \ge \epsilon 
 \end{split}
},
\end{equation}
where $D(y)$ is the matrix containing rows of $D$ corresponding to binding constraints at $y$. The point-to-point mapping $C:\reals^{mn} \times \reals^{mn} \to \reals^{mn}$ is defined such that $C(y, p) = y+ \gamma p$ where
\begin{equation}\label{eq:gammadef}
\gamma = \max\{ \gamma \in (0,1]: D(y+ \gamma p) \le d \}. 	
\end{equation}

The convergence theorem that we wish to establish is formally stated as follows.
\begin{lemma}\label{lm:conv}
	Suppose $U$ is concave and problem~\eqref{opt:abs} has a solution and its feasible set has a nonempty interior. Denote the optimal value of \eqref{opt:abs} by $U^\star$.  Then for any given feasible $y^0$,  any process $\{y^k\}_{k \in \Zbb_+}$ generated by algorithmic mapping $A_\epsilon$ in the sense of \eqref{eq:alg} has objective values $\{U(y^k)\}_{k \in \Zbb_+}$ such that 
	\begin{equation}
	U^\star - \lim_{k \to \infty} U(y^k)  \le \epsilon. 
	\end{equation}
\end{lemma}

Notice that in \eqref{eq:alg}, any point in the set $A_\epsilon(y^k)$ can be picked as $y^{k+1}$, and thus Lemma~\ref{lm:conv} is asserting a form of convergence for a family of an infinite numbers of processes $\{y^k\}_{k\in \Zbb_+}$. A classical result concerning with this type of convergence is Zangwill's \emph{global convergence} theorem \cite{zangwill1969nonlinear} (also see \cite{luenberger1984linear}):
\begin{theorem}
	Let $A$ be an algorithm on set $\mathcal X$, and suppose that starting from $x_0$ the sequence $\{x_k\}_{k \in \Zbb_+}$ is generated satisfying $x^{k+1} \in A(x^k)$. Let a solution set $\mathcal X^\star \subset \mathcal  X$ be given, and suppose
	\begin{enumerate}
		\item all points $x^k$ are contained in a compact set $\mathcal S \subset\mathcal X$, 
		\item there is a continuous function $U$ on $\mathcal X$ such that (i) if $x\not\in \mathcal X^\star$, then $U(y) > U(x)$ for all $y \in A(x)$, and (ii) if $x\in \mathcal X^\star$, then $U(y) \ge U(x)$ for all $y \in A(x)$,
		\item the mapping $A$ is closed at points outside of $\mathcal X^\star$. 
	\end{enumerate}
	Then the limit of any convergent subsequence of $\{x^k\}_{k \in \Zbb_+}$ is in $\mathcal X^\star$. 
\end{theorem}

The proof of Lemma~\ref{lm:conv} amounts to checking conditions in the theorem above for solution set $\Pcal^\star_\epsilon$ which are feasible points  $p$ for \eqref{opt:abs} such that $U^\star - U(p ) \le \epsilon$. 

\begin{IEEEproof}[Proof of Lemma~\ref{lm:conv}]
	The objective function $U$ is concave, hence continuous. The feasible set of \eqref{opt:abs} is an intersection of closed, compact, and convex sets and is also closed, compact and convex. It remains to show that the sequence $\{y^k\}_{k \in \Zbb_+}$ is feasible and ascent, and the mapping $A_\epsilon$ is closed. 
	
	\vspace{1mm}
\noindent{{\bf ($a$)} The sequence $\{y^k\}_{k\in \Zbb_+}$ is feasible for constraints of \eqref{opt:abs}.}
		
		The initial point $y^0$ is feasible by assumption. Suppose that $y^k$ is feasible. As $p^k \in G_\epsilon(y^k)$, we have $Ep^k=0$ and so $E(y^k  + \gamma p^k) =0$ for any $\gamma$. Therefore $Ey^{k+1}=0$. We also have $Dy^{k+1} \le d$ by the definition of $\gamma$ in \eqref{eq:gammadef}. (Notice that a nonzero $\gamma$ exists as the feasible set of \eqref{opt:abs} has nonempty interior.) 
		Finally, $y^{k+1}_i  = y^k_i + \gamma p^k_i \in \Pcal_i$ since $y_i^k \in \Pcal_i$, $y_i^k + p_i^k \in \Pcal_i$ and $\Pcal_i$ is convex. Thus $y^k$ is feasible for all $k$ by induciton. 
		
	\vspace{1mm}
\noindent{{\bf ($b$)}  The sequence $\{U(y^k)\}_{k\in \Zbb_+}$ is nondecreasing for $y^k \not \in \Pcal^\star_\epsilon$.} 

By the last condition in the definition of set $G_\epsilon$, we have $U(y^k + p^k) - U(y^k) \ge \epsilon >0$. As $\gamma \in (0,1]$ and $U$ is concave, we have $U(y^k + \gamma p^k) - U(y^k) > 0$.

	\vspace{1mm}
\noindent{{\bf ($c$)}   The mapping $A_\epsilon$ is closed.}

As the feasible set is closed, it suffices to prove that the graph of $A_\epsilon$ is closed. That is, for any sequence $\{(z^k, y^k)\}_{k \in \Zbb_+}$ with limit $(z,y)$ such that $z^k \in A_\epsilon (y^k)$, we need to show that $z\in A_\epsilon(y)$. It is easy to check that $z_i \in \Pcal_i$ as the set $\Pcal_i$ is closed for each $i$.
We thus proceed to show that the search direction defined by $p = z-y$ always satisfies linear constraints $Ep = 0$ and $D(y)p \le 0$. Indeed, the equality constraint holds by continuity of the linear mapping. Suppose that the inequality constraint $D(y) p \le 0$ does not hold. Then there exists a constraint, say the $\ell$th constraint, binding at $y$, i.e., $D_\ell^\top y = d_\ell$ such that $D_\ell^\top p = \delta >0$, where $D_\ell$ is the $\ell$th row of matrix $D$. By $y^k \to y$ and $D_\ell^\top y = d_\ell$, there exists a natural number $K_1$ such that for all $k \ge K_1$, 
\begin{equation}
D_\ell^\top y^k \ge d_\ell - \delta/4.
\end{equation}
Meanwhile, by $y^k \to y$ and $z^k \to z$, there exists a natural number $K_2$ such that for all $k \ge K_2$,
\begin{equation}
D_\ell^\top (z^k - y^k) \ge \delta/2 >0.
\end{equation}
It follows that for all $k \ge \max(K_1,K_2)$, we have 
\begin{equation}
D_\ell^\top z^k \ge D_\ell^\top y^k + \delta/2 \ge d_\ell + \delta/4 > d_\ell,
\end{equation}
contradicting to the fact $z^k$ is feasible as proved in item $(a)$, as $z^k \in A_\epsilon(y^k)$. 

It remains to show that there exists a $\eta \in [1, \infty)$ such that the search step before applying the mapping $C$ (curtailment), $\tilde p = \eta (z-y)$ is $\epsilon$-worthy. Suppose otherwise, then for all $\eta \in [1, \infty)$, 
\begin{equation}
U(y+ \eta (z-y)) - U(y ) \le \epsilon_1 < \epsilon.
\end{equation}
As $U$ is continuous, $y^k \to y$ and $z^k\to z$, there exists a $K_3$ such that for all $k \ge K_3$,
\begin{equation}
U(y^k + \eta (z^k-y^k)) - U(y^k) \le \epsilon_2, \quad \epsilon_2 \in (\epsilon_1, \epsilon),
\end{equation}
contradicting to $z^k \in A_\epsilon(y^k)$. Therefore $z\in A_\epsilon(y)$ and so $A_\epsilon$ is closed. 
\end{IEEEproof}

\begin{IEEEproof}[Proof of Theorem~\ref{thm:main}]
	We first recognize that \eqref{opt:swm} is a special case of \eqref{opt:abs}, with $n = S$, $m = |\Ical|$, $Ep = 0$ modeling the power balance constraints, $Dp \le d$ modeling the line capacity constraints. We further notice that the $\epsilon$-trading process (see footnote~\ref{ft:etrade}) is algorithm \eqref{eq:alg} with matrices $E$ and  $D$ suitably defined. Invoking Lemma~\ref{lm:conv} suggests that for any $\epsilon >0$, the trading state process $\{y^k\}_{k\in \Zbb_+}$ is such that $U^\star - \lim_{k \to \infty}U(y^k) \le \epsilon$ and therefore the claim in Theorem~\ref{thm:main} follows. 
\end{IEEEproof}

\begin{remark}[Trading process for distributed optimization]
	Given the general form of optimization~\eqref{opt:abs} and  per Lemma~\ref{lm:conv}, the trading process and its algorithmic correspondence~\eqref{eq:alg} define a framework for solving distributed optimization problems. Different from popular algorithms such as coordinate descent \cite{nesterov2012efficiency} and ADMM \cite{boyd2011distributed}, which iterate among coordinates with well-defined update order, the trading process would converge by updating some subsets of the coordinates (pairs of coordinates in the case of tree network or for network flow problems as discussed in Appendix~\ref{sec:tree}) according to any order,
 as long as conditions in the algorithmic mapping~\eqref{eq:alg} are met. 
	Another distinct feature of the trading process is that it does not specify a search direction for each update; it permits any search direction corresponding to a trade with suitable economic incentives ($\epsilon$-worthy trade) thus allows a great flexibility for designing platforms or systems in which agents could trade freely for their own benefits with an essential amount of coordination in place to ensure global constraints are satisfied.
\end{remark}

\section{Proof of Lemma~\ref{lemma:price}}
	Writing the Lagrangian of \eqref{opt:swm} and taking derivative with respect to $p_{i,s}$ for some RT participant $i$ gives
	\begin{equation}
	\rev{\P_i}(s) \frac{\partial \tilde u_{i,s}(p^\star)}{\partial p_{i,s}} + \lambda_{n,s}  - (\overline \eta_{i,s} - \underline \eta_{i,s}) = 0.
	\end{equation}
	Thus~\eqref{eq:lmp:gen} follows from this first order condition. When it is known that $p^\star_{i,s} \in \mathring{\Pcal_{i,s}}$, we have $\overline \eta_{i,s}^\star =  \underline \eta_{i,s}^\star =0$ by complementary slackness, and so equation~\eqref{eq:lmp:int} holds in this case.  


\section{Proof of Lemma 2}
Consider the optimization  of participant $i$ at bus $n$. The local constraint $p_i \in \Pcal_i$ can be expressed as
\begin{align}
& \underline p_i \le p_i \le \overline p_i, \notag\\
& p_i = \frac{1}{S}\ones \ones^\top p_i, \quad \mbox{if } i \in \IDA,
\end{align}
where the second constraint is the vector form of $p_{i,s} = \frac{1}{S} \sum_{s=1}^S p_{i,s}$, a convenient way to express the non-anticipation constraint. Let the dual variables for these constraints be denoted $\underline \eta_{i} \in \reals^S$, $\overline \eta_{i,s} \in \reals^S$ and $\zeta_{i}\in \reals^S$, respectively. Then the optimality condition for \eqref{opt:i} is  
\begin{subequations}\label{oc:i}
	\begin{align}
	& \lambda_{n,s} p_{i,s} +\rev{\P_i}(s)	\frac{\partial \tilde u_{i,s}(p_{i,s})}{\partial p_{i,s}} - (\overline \eta_{i,s} - \underline \eta_{i,s})  \nonumber \\
	& \quad \quad \quad \quad \quad \quad \quad \quad - (\ones_s - \frac{1}{S} \ones )^\top \zeta_i \mathbbm{1}_{i \in \IDA}= 0,\\
	 & \underline p_i \le p_i \le \overline p_i, \\
	 & p_i = \frac{1}{S} \ones \ones^\top p_i, \quad\quad \mbox{if } i \in \IDA,\\
	 & \underline \eta_i, \overline \eta_i \ge 0,
	\end{align}
\end{subequations}
where $\mathbbm{1}_{i\in \IDA}=1$ if $i \in \IDA$ and $0$ otherwise.  

Now consider the optimization for SO. Denote the dual variable for power balance constraint by $\gamma_s \in \reals$ and the dual variable for flow constraint by $\beta_s \in \reals^{L}$, $s\in [S]$. Then the optimality condition is 
\begin{subequations}\label{oc:so}
	\begin{align}
&		\lambda_{n,s} + \gamma_s + (H^\top \beta_s)_n =0, \quad n \in [N], s \in [S],\\
&		\ones^\top x_s = 0,\quad s \in [S], \\
&  Hx_s \le \overline f,    \,\,\quad  s\in [S],\\
& \beta_s \ge 0, \quad \,\,\,\quad s \in [S].
	\end{align}
\end{subequations}

Collecting optimality conditions~\eqref{oc:i} for all $i \in \Ical$ and that for SO \eqref{oc:so}, together with \eqref{eq:req}, we recover the optimality condition for \eqref{opt:swm}, where $\lambda_{n,s}$ is the dual for constraint. The claim in Lemma~\ref{lemma:eq} thus follows.

\section{Bilateral trading in tree network}\label{sec:tree}
In an example, \cite{Wu199975} demonstrates that multilateral trades involving more than two participants could be necessary when the network has cycles. The goal of this section is to complement that result by showing that when the network has no cycle,  bilateral trades are sufficient for the trading process to converge to the solution of centralized dispatch. 

For simplicity, we consider the deterministic case so that $S=1$ and $\Ical  = \IRT$. Without loss of generality, we further compress the notation by assuming there is only one participant connected to each node, so that $\Ical_n$ is a singleton for each $n$ and we use network index $n$ and participant index $i$ whichever is more convenient. These assumptions reduce our model to that of \cite{Wu199975}. 
The network is radial  so there is no cycle. With the DC approximation, the power flow model is equivalent to the standard network flow model \cite{bertsekas1998network} as shown in e.g. \cite{kim2016online}. 

We establish two decomposition results. We will state these results assuming that the accumulated trade of the network is zero ($x=0$) and we will work with the general line capacity constraints instead of constraints specified by the loading vector requirements. In a tree network, incorporating a nonzero accumulated trade $x$ amounts to modifying the line capacities and the utility functions; the loading vector requirements are special cases of the general line capacity constraints as they simply require the induced line flows on congested lines have nonpositive contribution to the congested direction.  Thus our treatment leads to no loss of generality. 

We first consider decomposing feasible multilateral trades into feasible bilateral trades.
\begin{proposition}\label{prop:f}
	Suppose that the line capacities are rational numbers, i.e., $\bar f \in \Qbb^L$. For any given feasible multilateral trade $p \in \Qbb^N$ involving more than two participants, i.e., $p \in \Pcal$ and $ \|p\|_0:=|\{n \in[N]: p_n \neq 0 \}| \ge 3$, there exists a finite number $K \in \Zbb_+$ of 
	bilateral trades $p^k \in \Qbb^N$ with $\|p^k\|_0=2$, $k=1, \dots, K$, that are \emph{sequentially feasible} such that $\sum_{k=1}^{K'} p^k \in \Pcal$ for any $K' \le K$ and satisfies $\sum_{k=1}^K p^k = p$.  
	
	Furthermore, under the same assumptions, there exists a finite number $\tilde K \in \Zbb_+$ of bilateral trades $\tilde p^k \in \Qbb^N$ with $\|\tilde p^k\|_0=2$, $k=1, \dots, \tilde K$, that are \emph{sequentially feasible for any ordering} satisfying $\sum_{k=1}^{\tilde K} p^k = p$. That is, let $\sigma: [\tilde K] \mapsto [\tilde K]$ be any permutation of the indices $1, \dots, \tilde K$, then $\sum_{k=1}^{\tilde{K}'} \tilde p^{\sigma(k)} \in \Pcal$, $\tilde K' < \tilde K$, and in particular $\tilde p^k \in \Pcal$ for any $k \in [\tilde K]$. 
\end{proposition}
\begin{IEEEproof}
	Let $(\Vcal, \Ecal)$ be the graph underlying the radial power network, with the node set $\Vcal = [N]$ and the edge set containing the $N-1$ lines of the network. Define the edge capacity to the capacity of the corresponding line. A flow on the graph is a $(N-1)$-vector $f$ that assigns a flow on each edge of the network satisfying flow conservation (for each node, in-flow equals the out-flow) and  edge capacity constraints.  Given a multilateral trade $p$, we denote the set of supply nodes by $\Vcal^+ = \{n \in \Vcal: p_n >0\}$ and the set of demand nodes by $\Vcal^- = \{n \in \Vcal:  p_n <0\}$. We then extend the graph by adding a source node $v^\mathrm{s}$ connecting to all the supply nodes and adding a sink node $v^\mathrm{t}$ connecting to all the demand node, so that the extended graph is 
	$(\tilde \Vcal, \tilde \Ecal)$ with $\tilde \Vcal = \Vcal \cup \{v^\mathrm{s}, v^\mathrm{t}\}$ and $\tilde \Ecal = \Ecal \cup \{(v^\mathrm{s}, v): v\in \Vcal^+\} \cup \{(v, v^\mathrm{t}): v\in \Vcal^-\}$. For edges $(v^\mathrm{s}, v)$, $v\in\Vcal^+$, we assign edge capacity to be $p_v>0$. Similarly for edges $(v, v^\mathrm{t})$, $v\in\Vcal^-$, the edge capacity is $- p_v >0$. 
	
	Now we make the observation that the multilateral trade $p$ and its induced power flow on the radial network is equivalent to the \emph{max flow} from $v^\mathrm{s}$ to $v^\mathrm{t}$ on the flow network $(\tilde \Vcal, \tilde \Ecal)$. In particular, the max flow solution would assign flows on the additional edges in $\tilde \Ecal \backslash \Ecal$ equal to the capacities. Together with the power flow induced by $p$, we obtain a feasible flow that maximizes the flow value from $v^\mathrm{s}$ to $v^\mathrm{t}$. Therefore, the problem of identifying feasible bilateral trades representing the multilateral trades is equivalent to finding feasible \emph{simple flows} representing the max flow on the flow network, where a simple flow from $v^\mathrm{s}$ to $v^\mathrm{t}$ is a flow on a simple path  from $v^\mathrm{s}$ to $v^\mathrm{t}$, which must contain exactly one supply node and one demand node. 
	
	The Ford-Fulkerson algorithm solves the max flow problem by iteratively identifying a feasible simple flow on the residual graph. By the definition  of the residual graph, this sequence of simple flows is sequentially feasible. Furthermore, it is known that Ford-Fulkerson terminates in a finite number of steps for rational inputs. We thus conclude that the finite collection of simple flows found by Ford-Fulkerson represents the finite collection of bilateral trades satisfying the requirements of the first part of the proposition. 
	
	 For the second part of the proposition, by Conformal Realization Theorem\cite[Proposition 1.1]{bertsekas1991linear} we know that there is a decomposition of the flow induced by $p$ on network $(\Vcal, \Ecal)$ into simple flows that are \emph{conformal} in the sense that the flow direction of the simple flows on each edge in $\Ecal$ is the same as that of the flow induced by $p$. 
	 It follows that this finite collection of simple flows is sequentially feasible for any ordering as for each edge if we choose the positive direction to be that of the flow direction induced by $p$, then for any $e\in \Ecal$, $\tilde f_e^k \ge 0$ and $\sum_{k=1}^{\tilde K} \tilde f_e = f_e$, where $f$ is the flow induced by $p$ and $\tilde f^k$ is the flow induced by $\tilde p^k$. The fact that each of the simple flows (and the corresponding bilateral trade) is feasible with respect to the original network constraint follows from choosing the $k$th simple flow as the first simple flow in the sequence. 
	 \end{IEEEproof}

We proceed to show that any \emph{non-redundant} profitable multilateral trades can be decomposed into a collection of profitable bilateral trades on a tree network.
\begin{definition}
	We say a  profitable  (and feasible) multilateral trade $p \in \reals^N$ contains redundancy if there exists a curtailment $\gamma \in [0,1]^N$ such that $\tilde p \in \reals^N$ defined by $\tilde p_n= \gamma_n p_n$, $n \in [N]$, is a feasible trade that achieves at least the same amount of profit as $p$. A profitable multilateral trade is deemed non-redundant if it does not contain redundancy. 
\end{definition}

Without loss of practicality, we focus  on and first state the result for the linear utility case so that $U_i(p_i) = \alpha_i p_i$ for some $\alpha_i \in \reals$, $i \in [N]$.  The extension to nonlinear case is discussed after that. 
\begin{proposition}\label{prop:p}
Under the same assumptions of Proposition~\ref{prop:f} and supposing that the utility function $U$ is linear,  any non-redundant profitable multilateral trades can be decomposed into a finite collection of profitable bilateral trades. Formally, given a non-redundant profitable trade $p \in \Qbb^n$, there exists a finite number $K$ of  bilateral trades $p^k$ that are sequentially feasible for any ordering 
and profitable, and satisfies $\sum_{k=1}^K p^k = p$. 
\end{proposition}
\begin{IEEEproof}
Consider the decomposition for the second part of Proposition~\ref{prop:f}. If all the bilateral trades in the decomposition are  profitable,  there is nothing to prove. Suppose there exists a bilateral trade  in the decomposition that is not profitable, denoted by $p'$. We claim that the remaining trade $p'' = p-p'$ is a feasible profitable trade that yields at least the same amount of profit as $p$ and therefore $p$ has redundancy. Indeed, by the proof of Proposition~\ref{prop:f}, $p''$ is feasible. Furthermore,
\begin{equation}
U(p'') - U(0) = \sum_{i = 1}^N \alpha_i p_i'' - 0 = \sum_{i=1}^N \alpha_i(p_i - p_i') > U(p) - U(0),
\end{equation}
as $U(p') = \sum_{i=1}^N \alpha_i < U(0) =0$ as $p'$ is not profitable.


\end{IEEEproof}

In general, when $U_i(p_i)$ is nonlinear but differentiable, we can decompose any given profitable trade $\bar p$ into $M\in\Zbb_+$ copies of trades $\bar p/M$, each of which is profitable by concavity of $U$.
 For $p_i \in [m\Delta  \bar p_i, (m+1) \Delta \bar p_i]$, $m = 0,\dots, M-1$, Taylor series offers a good linear approximation 
\begin{equation}
U_i( p_i) = U_i(m\Delta \bar p_i) + U'(m \Delta \bar p_i) p_i + o(\bar p_i/M), 
\end{equation}
where the last term denotes the approximation error which is of order higher than $\bar p_i/M$ and is negligible for practical purposes when $M$ is sufficiently large.  Applying Proposition~\ref{prop:p} gives a finite collection of profitable bilateral trades for each multilateral trade $\bar p/M$ with $m \bar p/M$ already scheduled into the system. Pooling these collections of bilateral trades gives a collection of bilateral trades that approximately represents the original multilateral trade. 

\rv{
\section{Trade verification and curtailment with local scenarios}\label{app:localscenario}
The bulk of the paper assumes that there is a commonly known set of \emph{global} scenarios $[S]$, which may not be practically available when participants are distributed over a large geographical area. In this section, we briefly discuss the other extreme setting where no global scenario is known a priori and participants are still allowed to submit contingent trades $p^k_s$, where $s\in \mathcal S^k$ with $\mathcal S^k$denoting the set of  \emph{local} scenarios  that is known to the participants involved in the $k$th trade.  Motivated by common sources of uncertainty in power systems (e.g. renewable generation level), we consider the setting that $\{p^k_s\in \reals^{|\Ical|}: s\in \mathcal S^k\}$ is an interval\footnote{\rv{When this does not hold, we can form the interval by finding the point-wise extremum of vector $p^k_s$.}} in $\reals^{|\Ical|}$. That is, the participants submit a lower bound $\underline p^k$  and an upper bound $\overline p^k$ to the SO, so that given any realization of the local uncertainty $s\in \mathcal S^k$, the resulting contingent trade satisfies $p^k_s \in [\underline p^k, \overline p^k]$. Since the SO may not be capable to identify the exact correlation among these local scenarios, the verification of power network constraints and curtailment has to be \emph{robust} with respect to any combinations of the local scenarios. That is, the accumulated network power injection defined by
\begin{equation}\label{eq:xkns}
	x^k_{n, s} = \sum_{\kappa=0}^{k-1} \sum_{i\in \Ical_n}\gamma^\kappa p^\kappa_{i,s}, \quad  n \in [N],
\end{equation}
must be feasible, i.e., 
\begin{equation}\label{eq:tocheck}
	x^k_s \in \mathcal P, 
\end{equation}
for
 every global scenario $s = (s^0, \dots, s^{k-1})$ generated by local scenarios  $s^\kappa \in \mathcal S^\kappa$, $\kappa = 0,\dots, k-1$. 
}

\rv{
We show that checking and curtailing a new trade  $p^k_s \in [\underline p^k, \overline p^k]$ with its corresponding network injection $q^k_s \in[\underline q^k, \overline q^k]$ can be done in an efficient way so that this process can be carried out inductively. Under our assumption that the set of contingent trades is an interval, the verification of $x^k_s + q^k_s \in \mathcal P$, for all $s = (s^0, \dots, s^k)$, $s^\kappa \in \mathcal S^\kappa$, is equivalent to checking
\begin{equation}
	H (x^k_s + q^k_s ) \le \bar f, 
\end{equation}
where $x^k_s$ is defined as in \eqref{eq:xkns} with $p^\kappa_s \in [\underline p^\kappa, \overline p^\kappa]$, $\kappa = 0, \dots, k-1$, and $q^k_s \in[\underline q^k, \overline q^k]$. 
Given the intervals for previous trades $p^\kappa_s \in [\underline p^\kappa, \overline p^\kappa]$, $\kappa = 0, \dots, k-1$, and their curtailment factors, we can form the corresponding interval $x^k_s \in [\underline x^k, \overline x^k]$. 
Then verifying the condition above  can be done by solving the following robust linear program:
\begin{subequations}
	\begin{align}
		\max_{\gamma \in [0,1]}\quad &\gamma \\
		\mbox{s.t.}\quad & x^k_{n, s} = \sum_{\kappa=0}^{k-1} \sum_{i\in \Ical_n}\gamma^\kappa p^\kappa_{i,s},  \quad n \in[N],\\
		&	H (x^k_s + \gamma q^k_s ) \le \bar f,
	\end{align}
\end{subequations}
where the last constraint must hold for all $x^k_s \in [\underline x^k, \overline x^k]$ and $q^k_s \in[\underline q^k, \overline q^k]$. In particular, if the optimal value of this program is $\gamma^\star =1$, then the new trade is feasible; if the optimal value is less than $1$, the new trade needs to be curtailed with curtailment factor $\gamma^k = \gamma^\star <1$. Finally, we note that this optimization can be solved efficiently by a bisection process (as the optimization variable is a scaler) equipped with a constraint checking sub-routine, which verifies the \emph{strong solvability of a set of interval linear inequalities} in polynomial time (cf. \cite[Section 2.13]{fiedler2006linear}).
}

\end{document}